\documentclass[times,sort&compress,3p]{elsarticle}
\usepackage[labelfont=bf]{caption}

\usepackage{amsmath,amsfonts,amssymb,amsthm,booktabs,color,epsfig,graphicx,url}
\usepackage{epstopdf}
\usepackage{diagbox}
\usepackage{slashbox}
\theoremstyle{plain}
\newtheorem{theorem}{Theorem}
\newtheorem{proposition}{Proposition}
\newtheorem{lemma}{Lemma}

\theoremstyle{definition}

\newtheorem{remark}{Remark}
\newtheorem{example}{Example}

\begin{document}

\begin{frontmatter}

\title{An analysis of multivariate measures of skewness and kurtosis of
skew-elliptical distributions}

\author[label1]{Baishuai  Zuo}
\author[label2]{Narayanaswamy Balakrishnan}
\author[label1]{Chuancun Yin\corref{mycorrespondingauthor}}
\address[label1]{School of Statistics and Data Science, Qufu Normal University, Qufu, Shandong 273165, P. R. China}
\address[label2]{Department of Mathematics and Statistics, McMaster University, Hamilton, Ontario, Canada}
\cortext[mycorrespondingauthor]{Corresponding author. Email address: \url{ccyin@qfnu.edu.cn (C. Yin)}}

\begin{abstract}
This paper examines eight measures of
skewness and Mardia measure of kurtosis for skew-elliptical distributions. Multivariate measures of
skewness considered include  Mardia, Malkovich-Afifi, Isogai, Song,
Balakrishnan-Brito-Quiroz, M$\acute{o}$ri, Rohatgi and Sz$\acute{e}$kely, Kollo and Srivastava measures. We first study
 the canonical form of skew-elliptical distributions, and then derive exact expressions of all measures of skewness and kurtosis for the family of skew-elliptical distributions, except for Song's measure. Specifically, the formulas of these measures for skew normal, skew $t$, skew logistic, skew Laplace, skew Pearson type II and skew Pearson type VII distributions are obtained. Next, as in Malkovich and Afifi (1973), test
statistics based on a random sample are constructed for illustrating the usefulness of the established results. In a Monte Carlo simulation study, different measures of
skewness and  kurtosis for $2$-dimensional skewed distributions are calculated and compared. Finally, real data is analyzed to demonstrate all the results.
\end{abstract}


\end{frontmatter}

\section{Introduction}
Mardia measures of skewness and kurtosis  for an arbitrary $k$-dimensional distribution $F$ with mean vector $\boldsymbol{\xi}$ and covariance matrix $\mathbf{\Delta}$ are defined, respectively, as
 \begin{align}\label{(a3)}
 \beta_{1,k}=\mathrm{E}[(\mathbf{X}-\boldsymbol{\xi})^{\top}\mathbf{\Delta}^{-1}(\mathbf{Y}-\boldsymbol{\xi})]^{3},
 \end{align}
 \begin{align}\label{(a5)}
 \beta_{2,k}=\mathrm{E}[(\mathbf{Y}-\boldsymbol{\xi})^{\top}\mathbf{\Delta}^{-1}(\mathbf{Y}-\boldsymbol{\xi})]^{2},
 \end{align}
 where $\mathbf{X}$ and $\mathbf{Y}$ are two independent and identically distributed random vectors from distribution $F$.
Since Mardia (1970) introduced these popular and commonly
used measures of multivariate skewness and kurtosis of an arbitrary $k$-dimensional distribution, many measures and indices of asymmetry have been proposed in the literature from different viewpoints. For example, Malkovich and Afifi (1973)  proposed a measure of multivariate skewness, as follows. Using the unit $k$-dimensional sphere $\phi_{k} = \{\boldsymbol{u} \in \mathbb{R}^{k}; ||\boldsymbol{u}|| = 1\},$ for $\boldsymbol{u} \in \phi_{k},$ the usual
univariate measure of skewness in the $\boldsymbol{u}$-direction is defined as
\begin{align*}
\beta_{1}(\boldsymbol{u})=\frac{[\mathrm{E}\{\boldsymbol{u}^{\top}(\mathbf{Y}-\mathrm{E}(\mathbf{Y}))\}^{3}]^{2}}{[\mathrm{Var}(\boldsymbol{u}^{\top}\mathbf{Y})]^{3}},
\end{align*}
and then a measure of skewness has been defined by them as
\begin{align*}
\beta_{1}^{\ast}=\sup_{\boldsymbol{u} \in \phi_{k}}\beta_{1}(\boldsymbol{u}).
\end{align*}

Based on  Mardia and Malkovich-Afifi measures, Balakrishnan et al. (2007) modified the Malkovich-Afifi measure to produce an overall vectorial measure of skewness as follows:
\begin{align}\label{(a15)}
T=\int _{\phi_{k}}\boldsymbol{u}c_{1}(\boldsymbol{u})\mathrm{d}\lambda(\boldsymbol{u}),
\end{align}
where $c_{1}(\boldsymbol{u}) =\mathrm{ E}[(\boldsymbol{u}^{\top}\mathbf{Z})^{3}]$, $\mathbf{Z}=\Delta^{-1/2}(\mathbf{Y}-\boldsymbol{\xi})$ with $\boldsymbol{\xi}=\mathrm{E}[\mathbf{Y}]$ and $\Delta=\mathrm{Var}(\mathbf{Y})$, and $\lambda$ is the rotationally invariant probability measure on the unit $k$-dimensional sphere
$\phi_{k} = \{\boldsymbol{u} \in \mathbb{R}^{k}; ||\boldsymbol{u}|| = 1\}.$
Isogai (1982) introduced an overall extension of Pearson measure of skewness for $\mathbf{Y}$ as
\begin{align}\label{(a12)}
 S_{I}=(\boldsymbol{\xi}-\mathbf{M_{0}})^{\top}h^{-1}(\mathbf{\Delta})(\boldsymbol{\xi}-\mathbf{M_{0}}),
 \end{align}
 where $\boldsymbol{\xi}$, $\mathbf{\Delta}$ and $\mathbf{M_{0}}$ are the mean, the covariance matrix  and  the mode of $\mathbf{Y}$, respectively. In addition, $h (\mathbf{\Delta})$ is an
``appropriate" function of the covariance matrix.
Srivastava (1984) defined  a measure of skewness for the multivariate vector $\mathbf{Y}$ based on the principal components $\mathbf{F} = \mathbf{\Upsilon} \mathbf{Y}$ as follows:
 \begin{align*}
 s_{1,k}^{2}=\frac{1}{k}\sum_{i=1}^{k}\left\{\frac{(F_{i}-\theta_{i})^{3}}{\lambda_{i}^{3/2}}\right\}^{2}=\frac{1}{k}\sum_{i=1}^{k}\left\{\frac{\mathrm{E}[\boldsymbol{\gamma}_{i}^{\top}(\mathbf{Y}-\boldsymbol{\xi})]^{3}}{\lambda_{i}^{3/2}}\right\}^{2},
 \end{align*}
where $\mathbf{\Upsilon}=(\boldsymbol{\gamma}_{1},\cdots,\boldsymbol{\gamma}_{k})$
denotes the matrix of eigenvectors of the covariance matrix $\Delta=\mathrm{Var}(\mathbf{Y})$, corresponding to the eigenvalues $\lambda_{1}, . . . , \lambda_{k}$, $F_{i}=\boldsymbol{\gamma}_{i}^{\top}\mathbf{Y}$ and $\theta_{i}=\boldsymbol{\gamma}_{i}^{\top}\boldsymbol{\xi}$ with $\boldsymbol{\xi}=\mathrm{E}[\mathbf{Y}]$.
M$\acute{o}$ri et al. (1994) proposed  a vectorial measure of skewness as a $k$-dimensional vector. For the standardized vector $\mathbf{Z}=\Delta^{-1/2}(\mathbf{Y}-\boldsymbol{\xi})$, this measure of skewness is defined as
 \begin{align*}
 s(\mathbf{Y})=\mathrm{E}[||\mathbf{Z}||\mathbf{Z}]=\sum_{i=1}^{k}\mathrm{E}[Z_{i}^{2}\mathbf{Z}]=\left(\sum_{i=1}^{k}\mathrm{E}[Z_{i}^{2}Z_{1}],\cdots,\sum_{i=1}^{k}\mathrm{E}[Z_{i}^{2}Z_{k}]\right)^{T}.
 \end{align*}
  Kollo (2008) noted that the M$\acute{o}$ri-Rohatgi-Sz$\acute{e}$kely measure of skewness dose not use all  the third-order mixed moments, and so he modified the M$\acute{o}$ri-Rohatgi-Sz$\acute{e}$kely measure of skewness to another measure of skewness including all mixed moments of third-order as follows:
  \begin{align*}
 b(\mathbf{Y})=\mathrm{E}\left[\sum_{i,j=1}^{k}Z_{i}Z_{j}\mathbf{Z}\right]=\left(\sum_{i,j=1}^{k}\mathrm{E}[Z_{i}Z_{j}Z_{1}],\cdots,\sum_{i,j=1}^{k}\mathrm{E}[Z_{i}Z_{j}Z_{k}]\right)^{T}.
 \end{align*}
 Using R$\acute{e}$nyi's entropy of order $\lambda$, Song (2001) introduced a measure of the shape of a distribution with pdf $f_{\mathbf{Y}}$ as
\begin{align*}
\mathfrak{L}_{R}(\lambda)=\frac{1}{1-\lambda}\log\int (f_{\mathbf{Y}}(\boldsymbol{y}))^{\lambda}\mathrm{d}\boldsymbol{y},
\end{align*}
and then a measure of skewness as
\begin{align*}
S(f_{\mathbf{Y}})=-2\frac{\mathrm{d}}{\mathrm{d}\lambda}\mathfrak{L}_{R}(1)=\mathrm{Var}[\log(f_{\mathbf{Y}}(\mathbf{Y}))].
\end{align*}
 Some other new distance based measure of asymmetry can be seen in the recent work of Baillien et al. (2022).

Balakrishnan and Scarpa (2012) mainly calculated and compared different measures of
multivariate skewness for the skew-normal family of distributions. Specifically, they derived the exact expressions of all the measures of
skewness, except for Song's measure, for the family of skew-normal distributions, and also presented
an approximate formula for Song's measure by using delta method. Kim and Zhao (2018) and Abdi et al. (2021) extended those results to scale mixtures
of skew-normal and mean-mixtures of multivariate normal distributions, respectively. Recently, Loperfido (2024) investigated skewness of multivariate mean-variance normal mixture distributions.
Yu and Yin (2023) derived first four moments of skew-elliptical distributions, and presented an approximate formula for Song's measure of skew-elliptical distribution, but is a complicated form. Balakrishnan et al. (2014) discussed a test of multivariate skew-normality by utilizing its canonical form.

Inspired by all these works, we calculate and compare here different measures of skewness and kurtosis for skew-elliptical distributions. Also,  we derive exact expressions for all the measures of
skewness and kurtosis for skew-elliptical distributions, except for Song's measure, but we use delta method and combined with some invariance properties to give an approximate formula for Song's measure for skew-elliptical distributions, which is a much simple form than that of Yu and Yin (2023).

The rest of the paper is organized as follows. Section 2 introduces the definitions and properties of skew-elliptical distributions. Section 3 presents the canonical form of skew-elliptical distributions, and gives an approximate formula for Song's measure and the exact expressions of other measures for skew-elliptical distributions. In Section 4, we consoder some special cases of elliptical distributions, are including skew normal, skew $t$, skew logistic, skew Laplace, skew Pearson type II and skew Pearson type VII distributions. In Section 5, we construct a test
statistic based on a random sample to illustrate the usefulness of all the established results. In Section 6, a simulation study is carried out. Analysis of real data is done in Sections 7. Finally, Section 8 presents some concluding remarks.
\section{Skew-elliptical distributions}
To define multivariate skew-elliptical distributions,
we first need to define multivariate elliptical distributions.  A $k$-dimensional random vector
$X = (X_{1}, \cdots , X_{k})^{\top}$ is said to have an elliptically symmetric distribution if its density function $f_{\mathbf{X}}(x)$ (when exists) has the form
\begin{align*}
f_{\mathbf{X}}(\boldsymbol{x})=|\mathbf{\Omega}|^{-\frac{1}{2}}
g^{(k)}\left((\boldsymbol{y-\mu})^{\top}\mathbf{\Omega}^{-1}(\boldsymbol{y-\mu})\right),
\end{align*}
where $\boldsymbol{\mu}\in\mathbb{R}^{k}$ is a location parameter, $\mathbf{\Omega}> 0$ is a $k\times k$ scale matrix, and $g^{(k)}$ is a spherical $k$-dimensional density generator; for details, see Fang et al. (1990).
A $k$-dimensional random vector $\mathbf{Y}$ is said to have a multivariate skew-elliptical distribution with a location parameter $\boldsymbol{\mu}\in\mathbb{R}^{k}$, a $k\times k$ scale matrix  $\mathbf{\Omega}> 0$ and shape parameter $\boldsymbol{\delta}\in\mathbb{R}^{k}$, denoted by $\mathbf{Y}\sim SE_{k} (\boldsymbol{\mu},\mathbf{\Omega},\boldsymbol{\delta}, g^{(k+1)})$, if its probability density function (pdf) is given  by (see Branco and Dey, 2001)
\begin{align*}
f_{\mathbf{Y}}(\boldsymbol{y})=2|\mathbf{\Omega}|^{-\frac{1}{2}}\int_{-\infty}^{\boldsymbol{\alpha}^{\top}(\boldsymbol{y-\mu})}
g^{(k+1)}\left(r^{2}+(\boldsymbol{y-\mu})^{\top}\mathbf{\Omega}^{-1}(\boldsymbol{y-\mu})\right)\mathrm{d}r,~\boldsymbol{y}\in\mathbb{R}^{k},
\end{align*}
where $\boldsymbol{\alpha}^{\top}=\frac{\boldsymbol{\delta}^{\top}\mathbf{\Omega}^{-1}}{(1-\boldsymbol{\delta}^{\top}\mathbf{\Omega}^{-1}\boldsymbol{\delta})^{1/2}}$, and $g^{(k+1)}$ is a spherical $k+1$-dimensional density generator.
 Assume $R$ is an absolutely continuous non-negative random variable, and $(U^{(1)},\mathbf{U}^{(p)\top})^{\top}\in\mathbb{R}^{k+1}$ has uniform distribution on the unit sphere and independent of $R$. Then, (Yu and Yin, 2023; Yin and  Balakrishnan, 2024)
 $$\mathbf{Y}\overset{d}{=}R(\boldsymbol{\delta}|U^{(1)}|+\Delta \mathbf{A}^{\top}\mathbf{U}^{(p)})+\boldsymbol{\mu},$$
  where $\mathbf{\Omega}=\mathbf{A}^{\top}\mathbf{A}$, $\Delta=\mathrm{diag}\left((1-\delta_{1}^{2})^{1/2},\cdots,(1-\delta_{k}^{2})^{1/2}\right)$, and $\delta_{i},~i=1,\cdots,k,$ denotes the $i$th component of $\boldsymbol{\delta}$.
 The density function of $R$ is
 $$h(r) =\frac{2\pi^{\frac{k+1}{2}}r^{k}g^{(k+1)}(r^{2})}{\Gamma(\frac{k+1}{2})},~r \in (0, +\infty),$$
 where $\Gamma(\cdot)$ denotes the Gamma function.
When $R=R_{0}\sim \chi_{k+1}$, and $\chi^{2}_{k+1}$ denotes a central chi square distribution with $k + 1$ degrees of freedom, $\mathbf{Y}$ will reduce to a skew-normal random vector (see, Genton, 2004).
\section{Main results}

\begin{lemma}\label{le.1}  Let $\mathbf{Y}\sim SE_{k} (\boldsymbol{\mu},\mathbf{\Omega},\boldsymbol{\delta}, g^{(k+1)})$,  $\mathbf{A}$ be a $k\times k$ non-singular matrix, and
$\boldsymbol{b}$ be a $k\times1$ vector. Then,
\begin{align}\label{(a1)}
\mathbf{AY}+\boldsymbol{b}\sim SE_{k} (\boldsymbol{\mu}_{\ast}, \mathbf{\Omega}_{\ast},\boldsymbol{\delta}_{\ast}, g^{(k+1)}),
\end{align}
where $\boldsymbol{\mu}_{\ast}=\mathbf{A}\boldsymbol{\mu}+\boldsymbol{b}$, $\mathbf{\Omega}_{\ast}=\mathbf{A}\mathbf{\Omega}\mathbf{A}^{\top}$ and $\boldsymbol{\delta}_{\ast}=\mathbf{A}\boldsymbol{\delta}$.
\end{lemma}
\noindent $\mathbf{Proof.}$ See Proposition 5.4 of Branco and Dey (2001). $\hfill\square$
\begin{remark}
In Proposition 5.4 of Branco and Dey (2001), one should have $\mathbf{\Omega}_{\ast}=\mathbf{A}\mathbf{\Omega}\mathbf{A}^{\top}$ instead of $\mathbf{\Omega}_{\ast}=\mathbf{A}^{\top}\mathbf{\Omega}\mathbf{A}$.
\end{remark}
With respect to canonical forms of distributions, Azzalini and  Capitanio (1999), Arellano-Valle et al. (2018), Capitanio (2020) and Abdi et al. (2021) presented canonical forms of skew-normal, scale and shape mixtures of multivariate skew-normal, scale mixtures of skew-normal and mean-mixtures of multivariate normal distributions, respectively.
Now, we derive the canonical form of skew-elliptical distributions.
\begin{theorem}\label{th.1}
Let $\mathbf{Y}\sim SE_{k} (\boldsymbol{\mu},\mathbf{\Omega},\boldsymbol{\delta}, g^{(k+1)})$. Then, there exists a linear transformation $\mathbf{Y}_{\ast}=\mathbf{A}_{\ast}(\mathbf{Y}-\boldsymbol{\mu})$
such that $\mathbf{Y}_{\ast}\sim SE_{k} (\boldsymbol{0},\mathbf{I}_{k},\boldsymbol{\delta}_{\ast}, g^{(k+1)})$, where at most one component of $\boldsymbol{\delta}_{\ast}$ is not zero, and $\boldsymbol{\delta}_{\ast} = (\delta_{\ast}, 0, \cdots , 0)^{\top}$ with $\delta_{\ast}=\boldsymbol{\delta}^{\top}\mathbf{\Omega}^{-1}\boldsymbol{\delta}$.
\end{theorem}
\noindent $\mathbf{Proof.}$ Because matrix $\mathbf{\Omega}$ is positive definite, there exists some invertible (non-singular) matrix $\mathbf{C}$ such that $\mathbf{\Omega}=\mathbf{C}^{\top}\mathbf{C}$. When $\boldsymbol{\delta}\neq\boldsymbol{0}$, there exists an orthogonal matrix $\mathbf{P}$ with the first column being proportional to $\mathbf{C}\mathbf{\Omega}^{-1}\boldsymbol{\delta}$, i.e. $\mathbf{P}=(\mathbf{C}\mathbf{\Omega}^{-1}\boldsymbol{\delta},P_{2},\cdots,P_{k})$. When $\boldsymbol{\delta}=\boldsymbol{0}$, we take $\mathbf{P}=\mathbf{I}_{k}$. Next, we let $\mathbf{A}_{\ast} =(\mathbf{C}^{-1}\mathbf{P})^{\top}$, and then use Lemma \ref{le.1} to obtain the desired result.  $\hfill\square$
\begin{remark}
In Theorem 5 of Abdi et al. (2021), there seems to be a small error; one should have   $\delta_{\ast}=\boldsymbol{\delta}^{\top}\mathbf{\Omega}^{-1}\boldsymbol{\delta}$ instead of $\delta_{\ast}=(\boldsymbol{\delta}^{\top}\mathbf{\Omega}^{-1}\boldsymbol{\delta})^{1/2}$. In addition, in (7) of  Azzalini and  Capitanio (1999),  $\alpha_{m}^{\ast} =(\boldsymbol{\alpha}^{\top}\mathbf{\Omega}\boldsymbol{\alpha})^{1/2}$ should be as $\alpha_{m}^{\ast} =\boldsymbol{\alpha}^{\top}\mathbf{\Omega}\boldsymbol{\alpha}$.
\end{remark}
\begin{theorem}\label{th.2}
Let $\mathbf{Y}_{\ast}\sim SE_{k} (\boldsymbol{0},\mathbf{I}_{k},\boldsymbol{\delta}_{\ast}, g^{(k+1)})$, where at most one component of $\boldsymbol{\delta}_{\ast}$ is not zero, and $\boldsymbol{\delta}_{\ast} = (\delta_{\ast}, 0, \cdots , 0)^{\top}$ with $\delta_{\ast}=\boldsymbol{\delta}^{\top}\mathbf{\Omega}^{-1}\boldsymbol{\delta}$. Then, the first four mixed moments of $\mathbf{Y}_{\ast}$ are given by
\begin{align*}
&\mathrm{E}[Y_{\ast 1}]=a\delta_{\ast},\\
&\mathrm{E}[Y_{\ast i}^{2}]=b,~i\in\{1,\cdots,k\},\\
&\mathrm{E}[Y_{\ast 1}^{3}]=c(3\delta_{\ast}-\delta_{\ast}^{3}),~E[Y_{\ast 1}Y_{\ast j}^{2}]=c\delta_{\ast},~j\in\{2,\cdots,k\},\\
&\mathrm{E}[Y_{\ast i}^{4}]=3d,~i\in\{1,\cdots,k\}, ~E[Y_{\ast i}^{2}Y_{\ast j}^{2}]=d,~i\neq j,~i,j\in\{1,\cdots,k\},\\
&\mathrm{E}[\cdot]=0,~otherwise,
\end{align*}
where
$a=\frac{\mathrm{E}[R]}{\mathrm{E}[R_{0}]}\sqrt{\frac{2}{\pi}}$, $b=\frac{\mathrm{E}[R^{2}]}{\mathrm{E}[R_{0}^{2}]}$, $c=\frac{\mathrm{E}[R^{3}]}{\mathrm{E}[R_{0}^{3}]}\sqrt{\frac{2}{\pi}}$ and $d=\frac{\mathrm{E}[R^{4}]}{\mathrm{E}[R_{0}^{4}]}$.
\end{theorem}
\noindent $\mathbf{Proof.}$ Using Theorem 3.1 of Yu and Yin (2023), after performing some algebra, we obtain the desired results. $\hfill\square$
\begin{remark}
As $R_{0}\sim \chi_{k+1}$, we readily obtain
$$\mathrm{E}[R_{0}]=\frac{\sqrt{2}\Gamma\left(\frac{k}{2}+1\right)}{\Gamma\left(\frac{k+1}{2}\right)},~\mathrm{E}[R_{0}^{2}]=k+1,~\mathrm{E}[R_{0}^{3}]=\frac{2\sqrt{2}\Gamma\left(\frac{k}{2}+2\right)}{\Gamma\left(\frac{k+1}{2}\right)},~\mathrm{E}[R_{0}^{4}]=(k+1)(k+3).$$
\end{remark}

\subsection{Mardia measures of skewness and kurtosis}

 \begin{proposition}\label{pro1}
Suppose $\mathbf{Y}\sim SE_{k} (\boldsymbol{\mu},\mathbf{\Omega},\boldsymbol{\delta}, g^{(k+1)})$, and the first three moments of $\mathbf{Y}$ exist. Then, we have
\begin{align}\label{(a4)}
\beta_{1,k}
 &=\frac{[3(c-ab)\delta_{\ast}+(2a^{3}-c)\delta_{\ast}^{3}]^{2}}{(b-a^{2}\delta_{\ast}^{2})^{3}}
 +\frac{3(k-1)(c-ab)^{2}\delta_{\ast}^{2}}{b^{2}(b-a^{2}\delta_{\ast}^{2})},
\end{align}
where  $\delta_{\ast}$, $a,~b ~ \mathrm{and}~c$ are the same as those in Theorem \ref{th.2}.
 \end{proposition}

\noindent $\mathbf{Proof.}$ Mardia (1970) pointed out that this measure of skewness is location and scale invariant (Mardia measure of skewness is invariant under nonsingular transformation $X=AY+b$).
 From Lemma 1 and Theorem 1, we know that the SE distribution is closed under affine transformations and has a canonical form. Now, let $\mathbf{Y}\sim SE_{k} (\boldsymbol{\mu},\mathbf{\Omega},\boldsymbol{\delta}, g^{(k+1)})$. Then, there exists a linear transformation $\mathbf{Y}_{\ast}=\mathbf{A}_{\ast}(\mathbf{Y}-\boldsymbol{\mu})$
such that $\mathbf{Y}_{\ast}\sim SE_{k} (\boldsymbol{0},\mathbf{I}_{k},\boldsymbol{\delta}_{\ast}, g^{(k+1)})$, where $\boldsymbol{\delta}_{\ast} = (\delta_{\ast}, 0, \cdots , 0)^{\top}$ with $\delta_{\ast}=\boldsymbol{\delta}^{\top}\mathbf{\Omega}^{-1}\boldsymbol{\delta}$. Without loss of generality, let $\mathbf{Y}_{\ast}=(Y_{\ast1},\cdots,Y_{\ast k})$.  We then compute Mardia measure of skewness of $\mathbf{Y}_{\ast}$.
  Let $\boldsymbol{\xi}_{\ast}=\mathrm{E}[\mathbf{X}_{\ast}]=\mathrm{E}[\mathbf{Y}_{\ast}]=\frac{\mathrm{E}[R]}{\mathrm{E}[R_{0}]}\sqrt{\frac{2}{\pi}}\boldsymbol{\delta}_{\ast} $ and $\Delta_{\ast}=\mathrm{Var}(\mathbf{X}_{\ast})=\mathrm{Var}(\mathbf{Y}_{\ast})=\frac{\mathrm{E}[R^{2}]}{\mathrm{E}[R_{0}^{2}]}\mathbf{I}_{k}-\frac{2\mathrm{E}^{2}[R]}{\pi\mathrm{E}^{2}[R_{0}]}\boldsymbol{\delta}_{\ast}\boldsymbol{\delta}_{\ast}^{\top}$ in (\ref{(a3)}). Then, the Mardia measure of skewness can be expressed as
 \begin{align*}
 \beta_{1,k}&=\mathrm{E}[(\mathbf{X}_{\ast}-\boldsymbol{\xi}_{\ast})^{\top}\mathbf{\Delta_{\ast}}^{-1}(\mathbf{Y}_{\ast}-\boldsymbol{\xi}_{\ast})]^{3}\\
 &=\frac{1}{(b-a^{2}\delta_{\ast}^{2})^{3}}\mathrm{E}^{2}[(Y_{\ast 1}-a\delta_{\ast})^{3}]+\frac{3}{b(b-a^{2}\delta_{\ast}^{2})^{2}}\sum_{i=2}^{k}\mathrm{E}^{2}[(Y_{\ast 1}-a\delta_{\ast})^{2}Y_{\ast i}]\\
 &~~~~+\frac{1}{b^{3}}\left\{\sum_{j=2}^{k}\sum_{i=2}^{k}\mathrm{E}^{2}[Y_{\ast j}Y_{\ast i}^{2}]+2\sum_{l=2}^{k}\sum_{2\leq i<j}^{k}\mathrm{E}^{2}[Y_{\ast i}Y_{\ast j}Y_{\ast l}] \right\}\\
 &~~~~+\frac{3}{b^{2}(b-a^{2}\delta_{\ast}^{2})}\left\{\sum_{i=2}^{k}\mathrm{E}^{2}[(Y_{\ast 1}-a\delta_{\ast})Y_{\ast i}^{2}]+2\sum_{2\leq i<j}^{k}\mathrm{E}^{2}[(Y_{\ast 1}-a\delta_{\ast})Y_{\ast i}Y_{\ast j}]     \right\}.
 \end{align*}
 From Theorem \ref{th.2}, we obtain
$$ \mathrm{E}^{2}[(Y_{\ast 1}-a\delta_{\ast})^{3}]=[3(c-ab)\delta_{\ast}+(2a^{3}-c)\delta_{\ast}^{3}]^{2},$$
$$\sum_{i=2}^{k}\mathrm{E}^{2}[(Y_{\ast 1}-a\delta_{\ast})^{2}Y_{\ast i}]=0,$$
$$\sum_{j=2}^{k}\sum_{i=2}^{k}\mathrm{E}^{2}[Y_{\ast j}Y_{\ast i}^{2}]+2\sum_{l=2}^{k}\sum_{2\leq i<j}^{k}\mathrm{E}^{2}[Y_{\ast i}Y_{\ast j}Y_{\ast l}] =0,$$
$$\sum_{i=2}^{k}\mathrm{E}^{2}[(Y_{\ast 1}-a\delta_{\ast})Y_{\ast i}^{2}]=(k-1)(c-ab)^{2}\delta_{\ast}^{2}$$
and
$$2\sum_{2\leq i<j}^{k}\mathrm{E}^{2}[(Y_{\ast 1}-a\delta_{\ast})Y_{\ast i}Y_{\ast j}]=0.$$
Hence, we obtain
\begin{align*}
\beta_{1,k}
 &=\frac{[3(c-ab)\delta_{\ast}+(2a^{3}-c)\delta_{\ast}^{3}]^{2}}{(b-a^{2}\delta_{\ast}^{2})^{3}}
 +\frac{3(k-1)(c-ab)^{2}\delta_{\ast}^{2}}{b^{2}(b-a^{2}\delta_{\ast}^{2})},
\end{align*}
completing the proof of the Proposition. $\hfill\square$
\begin{remark}
When $\boldsymbol{\delta}=\boldsymbol{0}$ in Proposition \ref{pro1}, $\mathbf{Y}$ will simply reduce to an elliptical vector. Then, its Mardia measure of skewness is $\beta_{1,k}=0$.
\end{remark}

\begin{proposition}\label{pro2}
Suppose $\mathbf{Y}\sim SE_{k} (\boldsymbol{\mu},\mathbf{\Omega},\boldsymbol{\delta}, g^{(k+1)})$, and  the first four moments of $\mathbf{Y}$ exist. Then, we have
\begin{align}\label{(a6)}
\beta_{2,k}
 &=\frac{3d+6a(ab-2c)\delta_{\ast}^{2}+a(4c-3a^{3})\delta_{\ast}^{4}}{(b-a^{2}\delta_{\ast}^{2})^{2}}
 +\frac{(k^{2}-1)d}{b^{2}}
 +\frac{2(k-1)[d+a(ab-2c)\delta_{\ast}^{2}]}{b(b-a^{2}\delta_{\ast}^{2})},
 \end{align}
where  $\delta_{\ast}$, $a,~b, ~c~ \mathrm{and}~d$ are the same as those in Theorem \ref{th.2}.
\end{proposition}
\noindent $\mathbf{Proof.}$ Mardia (1970) pointed out that
this measure of kurtosis is also location and scale invariant (Mardia measure of kurtosis is also invariant under nonsingular transformation $X=AY+b$). Let $\mathbf{Y}\sim SE_{k} (\boldsymbol{\mu},\mathbf{\Omega},\boldsymbol{\delta}, g^{(k+1)})$. From Lemma 1 and Theorem 1, we know that there exists a linear transformation $\mathbf{Y}_{\ast}=\mathbf{A}_{\ast}(\mathbf{Y}-\boldsymbol{\mu})$
such that $\mathbf{Y}_{\ast}\sim SE_{k} (\boldsymbol{0},\mathbf{I}_{k},\boldsymbol{\delta}_{\ast}, g^{(k+1)})$, where $\boldsymbol{\delta}_{\ast} = (\delta_{\ast}, 0, \cdots , 0)^{\top}$ with $\delta_{\ast}=\boldsymbol{\delta}^{\top}\mathbf{\Omega}^{-1}\boldsymbol{\delta}$. Without loss of generality, let $\mathbf{Y}_{\ast}=(Y_{\ast1},\cdots,Y_{\ast k})$. We then compute Mardia measure of kurtosis of $\mathbf{Y}_{\ast}$. For this purpose, let $\boldsymbol{\xi}_{\ast}=\mathrm{E}[\mathbf{Y}_{\ast}]=\frac{\mathrm{E}[R]}{\mathrm{E}[R_{0}]}\sqrt{\frac{2}{\pi}}\boldsymbol{\delta}_{\ast} $ and $\Delta_{\ast}=\mathrm{Var}(\mathbf{Y}_{\ast})=\frac{\mathrm{E}[R^{2}]}{\mathrm{E}[R_{0}^{2}]}\mathbf{I}_{k}-\frac{2\mathrm{E}^{2}[R]}{\pi\mathrm{E}^{2}[R_{0}]}\boldsymbol{\delta}_{\ast}\boldsymbol{\delta}_{\ast}^{\top}$ in (\ref{(a5)}). Then,
 the Mardia measure of kurtosis of $\mathbf{Y}_{\ast}$ can be expressed as
\begin{align*}
 \beta_{2,k}&=\mathrm{E}[(\mathbf{Y}_{\ast}-\boldsymbol{\xi}_{\ast})^{\top}\mathbf{\Delta}^{-1}(\mathbf{Y}_{\ast}-\boldsymbol{\xi}_{\ast})]^{2}\\
 &=\frac{1}{(b-a^{2}\delta_{\ast}^{2})^{2}}\mathrm{E}[(Y_{\ast 1}-a\delta_{\ast})^{4}]
 +\frac{1}{b^{2}}\left\{\sum_{i=2}^{k}\mathrm{E}[Y_{\ast i}^{4}]+2\sum_{2\leq i<j}^{k}\mathrm{E}[Y_{\ast i}^{2}Y_{\ast j}^{2}]\right\}
 +\frac{2}{b(b-a^{2}\delta_{\ast}^{2})}\sum_{i=2}^{k}\mathrm{E}[(Y_{\ast 1}-a\delta_{\ast})^{2}Y_{\ast i}^{2}]\\
 &=\frac{1}{(b-a^{2}\delta_{\ast}^{2})^{2}}\mathrm{E}[Y_{\ast 1}^{4}-4a\delta_{\ast}Y_{\ast 1}^{3}+6a^{2}\delta_{\ast}^{2}Y_{\ast 1}^{2}-4a^{3}\delta_{\ast}^{3}Y_{\ast 1}+a^{4}\delta_{\ast}^{4}]
 +\frac{1}{b^{2}}\left\{\sum_{i=2}^{k}\mathrm{E}[Y_{\ast i}^{4}]+2\sum_{2\leq i<j}^{k}\mathrm{E}[Y_{\ast i}^{2}Y_{\ast j}^{2}]\right\}\\
 &~~~~+\frac{2}{b(b-a^{2}\delta_{\ast}^{2})}\sum_{i=2}^{k}\mathrm{E}[Y_{\ast 1}^{2}Y_{\ast i}^{2}+a^{2}\delta_{\ast}^{2}Y_{\ast i}^{2}-2a\delta_{\ast}Y_{\ast 1}Y_{\ast i}^{2}].
 \end{align*}
 From Theorem \ref{th.2}, we obtain
 $$\mathrm{E}[Y_{\ast 1}^{4}-4a\delta_{\ast}Y_{\ast 1}^{3}+6a^{2}\delta_{\ast}^{2}Y_{\ast 1}^{2}-4a^{3}\delta_{\ast}^{3}Y_{\ast 1}+a^{4}\delta_{\ast}^{4}]=3d+6a(ab-2c)\delta_{\ast}^{2}+a(4c-3a^{3})\delta_{\ast}^{4},$$
 $$\sum_{i=2}^{k}\mathrm{E}[Y_{\ast i}^{4}]+2\sum_{2\leq i<j}^{k}\mathrm{E}[Y_{\ast i}^{2}Y_{\ast j}^{2}]=(k^{2}-1)d$$
and
 $$\sum_{i=2}^{k}\mathrm{E}[Y_{\ast 1}^{2}Y_{\ast i}^{2}+a^{2}\delta_{\ast}^{2}Y_{\ast i}^{2}-2a\delta_{\ast}Y_{\ast 1}Y_{\ast i}^{2}]=(k-1)[d+a(ab-2c)\delta_{\ast}^{2}].$$
 Hence, we obtain
 \begin{align*}
 \beta_{2,k}
 &=\frac{3d+6a(ab-2c)\delta_{\ast}^{2}+a(4c-3a^{3})\delta_{\ast}^{4}}{(b-a^{2}\delta_{\ast}^{2})^{2}}
 +\frac{(k^{2}-1)d}{b^{2}}
 +\frac{2(k-1)[d+a(ab-2c)\delta_{\ast}^{2}]}{b(b-a^{2}\delta_{\ast}^{2})},
 \end{align*}
 completing the proof of the Proposition. $\hfill\square$
 \begin{remark}
When $\boldsymbol{\delta}=\boldsymbol{0}$ in Proposition \ref{pro2}, $\mathbf{Y}$ will simply reduce to an elliptical vector. Then its Mardia measure of kurtosis reduces to
\begin{align}\label{(a7)}
\beta_{2,k}=\frac{(k+2)kd}{b^{2}}.
\end{align}
\end{remark}
In particular, Mardia's kurtosis measure for the multivariate $t$ distribution is given by
$$
\beta_{2,k}=\frac{(k+2)k(m-2)}{m-4},~ m>4,$$
which coincides with the expression in Eq. (16) of  Zografos (2008).

As Mardia's kurtosis measure of the normal distribution is $\beta_{2,k}=k(k+2)$, Mardia  (1974) proposed to measure excess kurtosis by the difference
\begin{align}\label{(a8)}
\gamma_{2,k}=\beta_{2,k}-k(k+2).
\end{align}

  M$\acute{o}$ri et al. (1994), Kollo and Srivastava (2005) and Ogasawara (2017) all  independently pointed out that Mardia's kurtosis is
never smaller than the sum of Mardia's skewness and the vector's dimension, i.e.,
\begin{align}\label{(a9)}
\beta_{2,k}\geq\beta_{1,k}+k.
\end{align}
This result is also verified  later in the simulation study.
\subsection{ Malkovich-Afifi measure of skewness}

\begin{proposition}\label{pro3}
Suppose $\mathbf{Y}\sim SE_{k} (\boldsymbol{\mu},\mathbf{\Omega},\boldsymbol{\delta}, g^{(k+1)})$, and  the first three moments of $\mathbf{Y}$ exist. Then, the Malkovich-Afifi measure of skewness is given by
\begin{align}\label{(a10)}
\beta_{1}^{\ast}&=\frac{\left[3(c-ab)\delta_{\ast}+(2a^{3}-c)\delta_{\ast}^{3}\right]^{2}}{\left[b-a^{2}\delta_{\ast}^{2}\right]^{3}},
\end{align}
where  $\delta_{\ast}$, $a,~b, ~c~ \mathrm{and}~d$ are the same as those in Theorem \ref{th.2}.
\end{proposition}
 \noindent $\mathbf{Proof.}$ Using Theorem 1, we obtain $\mathbf{Y}_{\ast}=(Y_{\ast 1}, \cdots,Y_{\ast k})^{\top}\sim SE_{k} (\boldsymbol{0},\mathbf{I}_{k},\boldsymbol{\delta}_{\ast}, g^{(k+1)})$, where $\boldsymbol{\delta}_{\ast} = (\delta_{\ast}, 0, \cdots , 0)^{\top}$ with $\delta_{\ast}=\boldsymbol{\delta}^{\top}\mathbf{\Omega}^{-1}\boldsymbol{\delta}$.
As Malkovich-Afifi measure of multivariate skewness is also location and scale invariant, we can consider its equivalent version as
\begin{align*}
\beta_{1}^{\ast}&=\sup_{\boldsymbol{u} \in \phi_{k}}\frac{[\mathrm{E}\{\boldsymbol{u}^{\top}(\mathbf{Y_{\ast}}-\mathrm{E}(\mathbf{Y_{\ast}}))\}^{3}]^{2}}{[\mathrm{Var}(\boldsymbol{u}^{\top}\mathbf{Y}_{\ast})]^{3}}.
\end{align*}
From Theorem \ref{th.2}, after performing some algebra, we obtain
\begin{align*}
\beta_{1}^{\ast}&=\sup_{\boldsymbol{u} \in \phi_{k}}\frac{\left[(2a^{3}-c)\delta_{\ast}^{3}u_{1}^{3}+3(c-ab)\delta_{\ast}u_{1}\right]^{2}}{\left[b-a^{2}\delta_{\ast}^{2}u_{1}^{2}\right]^{3}}\\
&=\frac{[\mathrm{E}\{Y_{\ast 1}-\mathrm{E}(Y_{\ast1})\}^{3}]^{2}}{[\mathrm{Var}(Y_{\ast 1})]^{3}}\\
&=\frac{\left[(2a^{3}-c)\delta_{\ast}^{3}+3(c-ab)\delta_{\ast}\right]^{2}}{\left[b-a^{2}\delta_{\ast}^{2}\right]^{3}},
\end{align*}
competing the proof of the Proposition. $\hfill\square$
\begin{remark}
When $\boldsymbol{\delta}=\boldsymbol{0}$ in Proposition \ref{pro3}, $\mathbf{Y}$ will simply reduce to an elliptical vector. Then, its  Malkovich-Afifi measure of skewness is $\beta_{1}^{\ast}=0$.
\end{remark}

\subsection{ Isogai measure of skewness}
  For obtaining this measure, we choose $h(\cdot)$ to be the identity function.
\begin{proposition}\label{pro5}
Suppose $\mathbf{Y}\sim SE_{k} (\boldsymbol{\mu},\mathbf{\Omega},\boldsymbol{\delta}, g^{(k+1)})$, and  the first three moments of $\mathbf{Y}$ exist. Then, the Isogai measure of skewness is given by
\begin{align}\label{(a13)}
 S_{I}&=\frac{(a\delta_{\ast}-m_{\ast 0})^{2}}{b-a^{2}\delta_{\ast}^{2}},
 \end{align}
 where $m_{\ast0}$ is the unique solution of the equation
\begin{align*}
2y_{1}\int_{-\infty}^{\eta(y_{1})}
g^{(k+1)'}\left(r^{2}+y_{1}^{2}\right)\mathrm{d}r+\eta(1)g^{(k+1)}\left(\eta^{2}(y_{1})+y_{1}^{2}\right)=0
\end{align*}
with $\eta(y_{1})=\frac{\delta_{\ast}y_{1}}{\sqrt{1-\delta_{\ast}^{2}}}$. In addition, $g^{(k+1)'}$ is the derivative of $g^{(k+1)}$, and $a$, $b$ and $\delta_{\ast}$ are the same as those in Theorem \ref{th.2}.
\end{proposition}
 \noindent $\mathbf{Proof.}$
 Using Theorem 1, we have $\mathbf{Y}_{\ast}=(Y_{\ast 1}, \cdots,Y_{\ast k})^{\top}\sim SE_{k} (\boldsymbol{0},\mathbf{I}_{k},\boldsymbol{\delta}_{\ast}, g^{(k+1)})$, where $\boldsymbol{\delta}_{\ast} = (\delta_{\ast}, 0, \cdots , 0)^{\top}$ with $\delta_{\ast}=\boldsymbol{\delta}^{\top}\mathbf{\Omega}^{-1}\boldsymbol{\delta}$. As this measure is location and
scale invariant,
let $\boldsymbol{\xi}_{\ast}=\mathrm{E}[\mathbf{Y}_{\ast}]=a\boldsymbol{\delta}_{\ast} $ and $\Delta_{\ast}=\mathrm{Var}(\mathbf{Y}_{\ast})=b\mathbf{I}_{k}-a^{2}\boldsymbol{\delta}_{\ast}\boldsymbol{\delta}_{\ast}^{\top}$ in (\ref{(a12)}).
 We then get
\begin{align*}
 S_{I}&=(\boldsymbol{\xi}_{\ast}-\mathbf{M_{\ast0}})^{\top}\mathbf{\Delta}_{\ast}^{-1}(\boldsymbol{\xi}_{\ast}-\mathbf{M_{\ast0}})\\
 &=\frac{(a\delta_{\ast}-M_{\ast 0 1})^{2}}{b-a^{2}\delta_{\ast}^{2}}+\frac{1}{b}\sum_{i=2}^{k}M_{\ast 0 i}^{2},
 \end{align*}
 where $\mathbf{M_{\ast0}}=(M_{\ast 0 1},\cdots,M_{\ast 0 k})^{\top}$ denotes the mode of $\mathbf{Y}_{\ast}$.\\
 As the pdf of $\mathbf{Y}_{\ast}$ is
  \begin{align*}
f_{\mathbf{Y}_{\ast}}(\boldsymbol{y})=2\int_{-\infty}^{\eta(y_{1})}
g^{(k+1)}\left(r^{2}+\boldsymbol{y}^{\top}\boldsymbol{y}\right)\mathrm{d}r,
\end{align*}
 we calculate the
mode $\mathbf{M_{\ast0}}$ by imposing the gradient of the density function to be equal to the null vector as follows:
\begin{align*}
&\frac{\partial f_{\mathbf{Y}_{\ast}}(\boldsymbol{y})}{\partial y_{1}}=4y_{1}\int_{-\infty}^{\eta(y_{1})}
g^{(k+1)'}\left(r^{2}+\boldsymbol{y}^{\top}\boldsymbol{y}\right)\mathrm{d}r+2\eta(1)g^{(k+1)}\left((\eta(y_{1}))^{2}+\boldsymbol{y}^{\top}\boldsymbol{y}\right)=0,\\
&\frac{\partial f_{\mathbf{Y}_{\ast}}(\boldsymbol{y})}{\partial y_{2}}=4y_{2}\int_{-\infty}^{\eta(y_{1})}
g^{(k+1)'}\left(r^{2}+\boldsymbol{y}^{\top}\boldsymbol{y}\right)\mathrm{d}r=0,\\
&\cdots\\
&\frac{\partial f_{\mathbf{Y}_{\ast}}(\boldsymbol{y})}{\partial y_{k}}=4y_{k}\int_{-\infty}^{\eta(y_{1})}
g^{(k+1)'}\left(r^{2}+\boldsymbol{y}^{\top}\boldsymbol{y}\right)\mathrm{d}r=0.
\end{align*}
 Hence, we obtain $\mathbf{M_{\ast0}}=(m_{\ast0},0,\cdots,0)^{\top}$, completing the proof of the Proposition. $\hfill\square$\\
 The vectorial measure, given similarly by Balakrishnan and Scarpa (2012), Kim and Zhao(2018) and Capitanio (2020), is
  \begin{align*}
 S_{C}&=\left(a-\frac{m_{\ast 0 }}{\delta_{\ast}}\right)\boldsymbol{\delta},
 \end{align*}
 where the direction of $\boldsymbol{\delta}$ can be
regarded as a measure of vectorial skewness of the SE distribution.
\begin{remark}
When $\boldsymbol{\delta}=\boldsymbol{0}$ in Proposition \ref{pro5}, $\mathbf{Y}$ will simply reduce to an elliptical vector. Then, its Isogai measure of skewness is given by
\begin{align}\label{(a14)}
 S_{I}&=0.
 \end{align}
\end{remark}
\subsection{Song measure}
It is difficult to derive an exact expression of Song measure for skew-elliptical distribution. So, we now present an approximate formula for this measure for the family of skew-elliptical distributions.
\begin{proposition}\label{pro6}
Suppose $\mathbf{Y}\sim SE_{k} (\boldsymbol{\mu},\mathbf{\Omega},\boldsymbol{\delta}, g^{(k+1)})$, and  the first four moments of $\mathbf{Y}$ exist. Then, an approximation of Song measure is given by
\begin{align}\label{(a19)}
S(f_{\mathbf{Y}})\approx(b-a^{2}\delta_{\ast}^{2})h_{\ast}^{2},
\end{align}
where
\begin{align*}
h_{\ast}=\frac{2a\delta_{\ast}\int_{-\infty}^{\eta(a\delta_{\ast})}
g^{(k+1)'}\left(r^{2}+a^{2}\delta_{\ast}^{2}\right)\mathrm{d}r+\eta(1)g^{(k+1)}\left((\eta(a\delta_{\ast}))^{2}+a^{2}\delta_{\ast}^{2}\right)}{\int_{-\infty}^{\eta(a\delta_{\ast})}
g^{(k+1)}\left(r^{2}+a^{2}\delta_{\ast}^{2}\right)\mathrm{d}r}
\end{align*}
with $\eta(x)=\frac{\delta_{\ast}x}{\sqrt{1-\delta_{\ast}^{2}}}$. In addition,
$a$, $b$ and $\delta_{\ast}$ are the same as those in Theorem \ref{th.2}.
\end{proposition}
\noindent $\mathbf{Proof.}$ To obtain  an approximation,  using the delta method, we have
\begin{align*}
\mathrm{Var}[\log(f_{\mathbf{Y}}(\mathbf{Y}))]=\mathrm{Var}[G_{\mathbf{Y}}(\mathbf{Y})]\approx [\nabla G_{\mathbf{Y}}(\mathrm{E}(\mathbf{Y}))]^{\top}\mathrm{Var}(\mathbf{Y})[\nabla G_{\mathbf{Y}}(\mathrm{E}(\mathbf{Y}))].
\end{align*}
As Song measure is also location and scale invariant, using Theorem \ref{th.1}, we have
 \begin{align}\label{(a17)}
S(f_{\mathbf{Y}})&=S(f_{\mathbf{Y_{\ast}}})\approx [\nabla G_{\mathbf{Y_{\ast}}}(\mathrm{E}(\mathbf{Y_{\ast}}))]^{\top}\mathrm{Var}(\mathbf{Y_{\ast}})[\nabla G_{\mathbf{Y_{\ast}}}(\mathrm{E}(\mathbf{Y_{\ast}}))],
\end{align}
where $\mathrm{E}[\mathbf{Y_{\ast}}]=a\boldsymbol{\delta}_{\ast}$ and $\mathrm{Var}(\mathbf{Y_{\ast}})=b\mathbf{I}_{k}-a^{2}\boldsymbol{\delta}_{\ast}\boldsymbol{\delta}_{\ast}^{\top}$.\\
While $G_{\mathbf{Y_{\ast}}}(\boldsymbol{y})=\log2+\log\int_{-\infty}^{\eta(y_{1})}
g^{(k+1)}\left(r^{2}+\boldsymbol{y}^{\top}\boldsymbol{y}\right)\mathrm{d}r$, we get
 \begin{align}\label{(a18)}
 \nabla G_{\mathbf{Y_{\ast}}}(\boldsymbol{y})=\left(\frac{\partial G_{\mathbf{Y_{\ast}}}(\boldsymbol{y})}{\partial y_{1}},\cdots,\frac{\partial G_{\mathbf{Y_{\ast}}}(\boldsymbol{y})}{\partial y_{k}}\right)^{\top},
  \end{align}
where
\begin{align*}
&\frac{\partial G_{\mathbf{Y_{\ast}}}(\boldsymbol{y})}{\partial y_{1}}=\frac{2y_{1}\int_{-\infty}^{\eta(y_{1})}
g^{(k+1)'}\left(r^{2}+\boldsymbol{y}^{\top}\boldsymbol{y}\right)\mathrm{d}r+\eta(1)g^{(k+1)}\left((\eta(y_{1}))^{2}+\boldsymbol{y}^{\top}\boldsymbol{y}\right)}{\int_{-\infty}^{\eta(y_{1})}
g^{(k+1)}\left(r^{2}+\boldsymbol{y}^{\top}\boldsymbol{y}\right)\mathrm{d}r},\\
&\frac{\partial G_{\mathbf{Y_{\ast}}}(\boldsymbol{y})}{\partial y_{2}}=\frac{2y_{2}\int_{-\infty}^{\eta(y_{1})}
g^{(k+1)'}\left(r^{2}+\boldsymbol{y}^{\top}\boldsymbol{y}\right)\mathrm{d}r}{\int_{-\infty}^{\eta(y_{1})}
g^{(k+1)}\left(r^{2}+\boldsymbol{y}^{\top}\boldsymbol{y}\right)\mathrm{d}r},\\
&\cdots\\
&\frac{\partial G_{\mathbf{Y_{\ast}}}(\boldsymbol{y})}{\partial y_{k}}=\frac{2y_{k}\int_{-\infty}^{\eta(y_{1})}
g^{(k+1)'}\left(r^{2}+\boldsymbol{y}^{\top}\boldsymbol{y}\right)\mathrm{d}r}{\int_{-\infty}^{\eta(y_{1})}
g^{(k+1)}\left(r^{2}+\boldsymbol{y}^{\top}\boldsymbol{y}\right)\mathrm{d}r}.
\end{align*}
  Let $\boldsymbol{y}=\mathrm{E}[\mathbf{Y_{\ast}}]$ in (\ref{(a18)}). Then, by combining it with (\ref{(a17)}), and then performing some algebra, we readily obtain (\ref{(a19)}), completing the proof of the Proposition. $\hfill\square$
\begin{remark}
When $\boldsymbol{\delta}=\boldsymbol{0}$ in Proposition \ref{pro6}, $\mathbf{Y}$ will simply reduce to an elliptical vector. Then, its approximation of Song measure is given by
\begin{align*}
S(f_{\mathbf{Y}})\approx 0.
\end{align*}
\end{remark}

\subsection{ Balakrishnan-Brito-Quiroz measure of skewness}
 To obtain a single measure, Balakrishnan et al. (2007) proposed the quantity $Q = T^{\top} \Sigma_{T}^{-1}T $, where T is as in (\ref{(a15)}) and $\Sigma_{T}$ is the covariance matrix of $T$. As it is difficult to compute $\Sigma_{T}$, they use $\Sigma_{Z}$ in place of $\Sigma_{T}$ to obtain,  we obtain
 \begin{align}
 Q_{\ast} = T^{\top} \Sigma_{Z}^{-1}T=T^{\top}T ,
 \end{align} also as done in Balakrishnan and Scarpa (2012). From Abdi et al. (2021), we can obtain the elements of $T$ as follows:
 \begin{align}\label{(a16)}
 T_{r}=\frac{3}{k(k+2)}\mathrm{E}[Z_{r}^{3}]+3\sum_{i\neq r}\frac{1}{k(k+2)}\mathrm{E}[Z_{i}^{2}Z_{r}]=\frac{3}{k(k+2)}\sum_{i= 1}^{k}\mathrm{E}[Z_{i}^{2}Z_{r}],
 \end{align}
 where the required $\mathrm{E}[Z_{r}^{3}]$ and  $\mathrm{E}[Z_{i}^{2}Z_{r}]$ can be obtained from the  following theorem.
 \begin{theorem}\label{th.4}
Let $\mathbf{Y}\sim SE_{k} (\boldsymbol{\mu},\mathbf{\Omega},\boldsymbol{\delta}, g^{(k+1)})$, $\boldsymbol{\xi}=\mathrm{E}[\mathbf{Y}]$ and $\Delta=\mathrm{Var}(\mathbf{Y})$. Then, $\mathbf{Z}=\Delta^{-1/2}(\mathbf{Y}-\boldsymbol{\xi})\sim SE_{k} (-a\boldsymbol{\delta}_{Z},\mathbf{\Omega}_{Z},\boldsymbol{\delta}_{Z}, g^{(k+1)}).$\\
 Further, we have
\begin{align*}
&\Sigma_{Z}=\mathrm{Var}(\mathbf{Z})=\mathbf{I}_{k},\\
&\mathrm{E}[Z_{i}^{3}]=M_{3}[(i-1)k+i,i],i\in\{1,\cdots,k\},\\
&\mathrm{E}[Z_{i}^{2}Z_{j}]=M_{3}[(i-1)k+i,j],~i,j\in\{1,\cdots,k\},\\
&\mathrm{E}[Z_{i}Z_{j}Z_{r}]=M_{3}[(i-1)k+r,j],~i,j,r\in\{1,\cdots,k\},
\end{align*}
where
\begin{align*}
 M_{3}=&2a^{3}\boldsymbol{\delta}_{Z}\otimes\boldsymbol{\delta}_{Z}^{\top}\otimes\boldsymbol{\delta}_{Z}
 -ab[\mathbf{\Omega}_{Z}\otimes\boldsymbol{\delta}_{Z}+\boldsymbol{\delta}_{Z}\otimes\mathbf{\Omega}_{Z}+\mathrm{vec}(\mathbf{\Omega}_{Z})\otimes\boldsymbol{\delta}_{Z}^{\top}]\\
&+c[\boldsymbol{\delta}_{Z}\otimes\mathbf{\Omega}_{Z}+\mathrm{vec}(\mathbf{\Omega}_{Z})\boldsymbol{\delta}_{Z}^{\top}+(\mathbf{I}_{k}\otimes\boldsymbol{\delta}_{Z})\mathbf{\Omega}_{Z}-(\mathbf{I}_{k}\otimes\boldsymbol{\delta}_{Z})(\boldsymbol{\delta}_{Z}\otimes\boldsymbol{\delta}_{Z}^{\top})]
 \end{align*}
 with
 $$\mathbf{\Omega}_{Z}=(b\mathbf{\Omega}-a^{2}\boldsymbol{\delta}\boldsymbol{\delta}^{\top})^{-1/2}\mathbf{\Omega}(b\mathbf{\Omega}-a^{2}\boldsymbol{\delta}\boldsymbol{\delta}^{\top})^{-1/2},
 ~\boldsymbol{\delta}_{Z}=(b\mathbf{\Omega}-a^{2}\boldsymbol{\delta}\boldsymbol{\delta}^{\top})^{-1/2}\boldsymbol{\delta},$$
and $M_{3}[i,j]$ denotes the $i$th row and $j$th column element of $M_{3}$. In addition, $\otimes$ and $\mathrm{vec}$ are the Kronecker product and vectorizing operator, respectively (for details, see Fang and Zhang, 1990).
\end{theorem}
\noindent $\mathbf{Proof.}$ Using Theorem 5 of Yin and  Balakrishnan (2024), and after performing some algebra, we have
 $\boldsymbol{\xi}=\mathrm{E}[\mathbf{Y}]=\boldsymbol{\mu}+a\boldsymbol{\delta}$ and $\Delta=\mathrm{Var}(\mathbf{Y})=b\mathbf{\Omega}-a^{2}\boldsymbol{\delta}\boldsymbol{\delta}^{\top}$.
 Next, using Lemma \ref{le.1} and after some algebra, we obtain $\mathbf{Z}=\Delta^{-1/2}(\mathbf{Y}-\boldsymbol{\xi})\sim SE_{k} (\boldsymbol{\mu}_{Z},\mathbf{\Omega}_{Z},\boldsymbol{\delta}_{Z}, g^{(k+1)}),$
 with
 $$\boldsymbol{\mu}_{Z}=\Delta^{-1/2}(\boldsymbol{\mu}-\boldsymbol{\xi})=-a(b\mathbf{\Omega}-a^{2}\boldsymbol{\delta}\boldsymbol{\delta}^{\top})^{-1/2}\boldsymbol{\delta}=-a\boldsymbol{\delta}_{Z},$$
$$\mathbf{\Omega}_{Z}=\Delta^{-1/2}\mathbf{\Omega}\Delta^{-1/2}=(b\mathbf{\Omega}-a^{2}\boldsymbol{\delta}\boldsymbol{\delta}^{\top})^{-1/2}\mathbf{\Omega}(b\mathbf{\Omega}-a^{2}\boldsymbol{\delta}\boldsymbol{\delta}^{\top})^{-1/2},$$
$$\boldsymbol{\delta}_{Z}=\Delta^{-1/2}\boldsymbol{\delta}=(b\mathbf{\Omega}-a^{2}\boldsymbol{\delta}\boldsymbol{\delta}^{\top})^{-1/2}\boldsymbol{\delta}.$$
 Moreover, by Theorem 5 of Yin and  Balakrishnan (2024) again, we obtain
 \begin{align*}
 M_{1}&=\boldsymbol{\mu}_{Z}+a\boldsymbol{\delta}_{Z}=\boldsymbol{0},\\
 M_{2}&=\boldsymbol{\mu}_{Z}\boldsymbol{\mu}_{Z}^{T}+a(\boldsymbol{\mu}_{Z}\boldsymbol{\delta}_{Z}^{T}+\boldsymbol{\delta}\boldsymbol{\mu}_{Z}^{T})+b\mathbf{\Omega}_{Z}=\mathbf{I}_{k},
 \\
 \Sigma_{Z}&=\mathrm{Var}(\mathbf{Z})=M_{2}-M_{1}M_{1}^{T}=\mathbf{I}_{k},
 \\
 M_{3}=&\boldsymbol{\mu}_{Z}\otimes\boldsymbol{\mu}_{Z}^{\top}\otimes\boldsymbol{\mu}_{Z}
 +a[\boldsymbol{\delta}_{Z}\otimes\boldsymbol{\mu}_{Z}^{\top}\otimes\boldsymbol{\mu}_{Z}+\boldsymbol{\mu}_{Z}\otimes\boldsymbol{\delta}_{Z}^{\top}\otimes\boldsymbol{\mu}_{Z}+\boldsymbol{\mu}_{Z}\otimes\boldsymbol{\mu}_{Z}^{\top}\otimes\boldsymbol{\delta}_{Z}]\\
 &+b[\mathbf{\Omega}_{Z}\otimes\boldsymbol{\mu}_{Z}+\boldsymbol{\mu}_{Z}\otimes\mathbf{\Omega}_{Z}+vec(\mathbf{\Omega}_{Z})\otimes\boldsymbol{\mu}_{Z}^{\top}]\\
&+c[\boldsymbol{\delta}_{Z}\otimes\mathbf{\Omega}_{Z}+\mathrm{vec}(\mathbf{\Omega}_{Z})\boldsymbol{\delta}_{Z}^{\top}+(\mathbf{I}_{k}\otimes\boldsymbol{\delta}_{Z})\mathbf{\Omega}_{Z}-(\mathbf{I}_{k}\otimes\boldsymbol{\delta}_{Z})(\boldsymbol{\delta}_{Z}\otimes\boldsymbol{\delta}_{Z}^{\top})]\\
=&2a^{3}\boldsymbol{\delta}_{Z}\otimes\boldsymbol{\delta}_{Z}^{\top}\otimes\boldsymbol{\delta}_{Z}
 -ab[\mathbf{\Omega}_{Z}\otimes\boldsymbol{\delta}_{Z}+\boldsymbol{\delta}_{Z}\otimes\mathbf{\Omega}_{Z}+\mathrm{vec}(\mathbf{\Omega}_{Z})\otimes\boldsymbol{\delta}_{Z}^{\top}]\\
&+c[\boldsymbol{\delta}_{Z}\otimes\mathbf{\Omega}_{Z}+\mathrm{vec}(\mathbf{\Omega}_{Z})\boldsymbol{\delta}_{Z}^{\top}+(\mathbf{I}_{k}\otimes\boldsymbol{\delta}_{Z})\mathbf{\Omega}_{Z}-(\mathbf{I}_{k}\otimes\boldsymbol{\delta}_{Z})(\boldsymbol{\delta}_{Z}\otimes\boldsymbol{\delta}_{Z}^{\top})],
 \end{align*}
 completing the proof of the Theorem. $\hfill\square$

 \subsection{M$\acute{o}$ri-Rohatgi-Sz$\acute{e}$kely measure of skewness}
 M$\acute{o}$ri et al. (1994) proposed  a vectorial measure of skewness as a $k$-dimensional vector. For the standardized vector $\mathbf{Z}=\Delta^{-1/2}(\mathbf{Y}-\boldsymbol{\xi})$, this measure of skewness is defined as:
 \begin{align*}
 s(\mathbf{Y})=\mathrm{E}[||\mathbf{Z}||\mathbf{Z}]=\sum_{i=1}^{k}\mathrm{E}[Z_{i}^{2}\mathbf{Z}]=\left(\sum_{i=1}^{k}\mathrm{E}[Z_{i}^{2}Z_{1}],\cdots,\sum_{i=1}^{k}\mathrm{E}[Z_{i}^{2}Z_{k}]\right)^{\top}.
 \end{align*}
 From Theorem \ref{th.4}, we can calculate $\mathrm{E}[Z_{i}^{3}]$, $i\in\{1,\cdots,k\}$, and $\mathrm{E}[Z_{i}^{2}Z_{j}],$ $i\neq j$, $i,j\in\{1,\cdots,k\}$. Hence, $s(\mathbf{Y})$ can be readily found.
 \subsection{Kollo measure of skewness}
  Kollo (2008) noted that the M$\acute{o}$ri-Rohatgi-Sz$\acute{e}$kely measure of skewness does not use all  third-order mixed moments, and he therefore modified that measure of skewness to a measure of skewness by including all mixed moments of third-order as
  \begin{align*}
 b(\mathbf{Y})=\mathrm{E}\left[\sum_{i,j=1}^{k}Z_{i}Z_{j}\mathbf{Z}\right]=\left(\sum_{i,j=1}^{k}\mathrm{E}[Z_{i}Z_{j}Z_{1}],\cdots,\sum_{i,j=1}^{k}\mathrm{E}[Z_{i}Z_{j}Z_{k}]\right)^{\top}.
 \end{align*}
 From Theorem \ref{th.4}, we can compute $\mathrm{E}[Z_{i}^{3}]$, $i\in\{1,\cdots,k\}$,  $\mathrm{E}[Z_{i}^{2}Z_{j}],$ $i\neq j$, $i,j\in\{1,\cdots,k\}$ and $\mathrm{E}[Z_{i}Z_{j}Z_{l}],$ $i\neq j\neq l$, $i,j,l\in\{1,\cdots,k\}$, and so we can readily obtain $b(\mathbf{Y})$.
\subsection{Srivastava measure of skewness}
It is known that for any random vector $\mathbf{X}$, there is a relationship between central moments of third order ($\overline{M}_{3}(\mathbf{X})$) and  first three non-central moments ($\mathrm{E}(\mathbf{X})$, $M_{2}(\mathbf{X})$ and $M_{3}(\mathbf{X})$) as follows (see, for example, Kollo and Srivastava, 2005):
\begin{align*}
\overline{M}_{3}(\mathbf{X})=M_{3}(\mathbf{X})-M_{2}(\mathbf{X})\otimes\mathrm{E}(\mathbf{X})-\mathrm{E}(\mathbf{X})\otimes M_{2}(\mathbf{X})-\mathrm{vec}(M_{2}(\mathbf{X}))\mathrm{E}(\mathbf{X})^{\top}+2\mathrm{E}(\mathbf{X})\mathrm{E}(\mathbf{X})^{\top}\otimes\mathrm{E}(\mathbf{X}).
\end{align*}
From Theorem 5 of Yin and  Balakrishnan (2024), by using the
relations for affine transformations of moments (see Balakrishnan and Scarpa, 2012), we find
\begin{align*}
\mathrm{E}(F_{i})=\boldsymbol{\gamma}_{i}^{\top}\mathrm{E}[\mathbf{Y}],~~
M_{2}(F_{i})=\boldsymbol{\gamma}_{i}^{\top}\mathrm{E}[\mathbf{Y}\mathbf{Y}^{\top}]\boldsymbol{\gamma}_{i},~~
M_{3}(F_{i})=(\boldsymbol{\gamma}_{i}^{\top}\otimes\boldsymbol{\gamma}_{i}^{\top})M_{3}(\mathbf{Y})\boldsymbol{\gamma}_{i},
\end{align*}
and consequently
\begin{align*}
&\mathrm{E}[\boldsymbol{\gamma}_{i}^{\top}(\mathbf{Y}-\boldsymbol{\xi})]^{3}=\overline{M}_{3}(F_{i})\\
&=\left(\boldsymbol{\gamma}_{i}^{\top}\otimes\boldsymbol{\gamma}_{i}^{\top}\right)M_{3}(\mathbf{Y})\boldsymbol{\gamma}_{i}
-\left(\boldsymbol{\gamma}_{i}^{\top}\mathrm{E}[\mathbf{Y}\mathbf{Y}^{\top}]\boldsymbol{\gamma}_{i}\right)\otimes\left(\boldsymbol{\gamma}_{i}^{\top}\mathrm{E}[\mathbf{Y}]\right)
-\left(\boldsymbol{\gamma}_{i}^{\top}\mathrm{E}[\mathbf{Y}]\right)\otimes \left(\boldsymbol{\gamma}_{i}^{\top}\mathrm{E}[\mathbf{Y}\mathbf{Y}^{\top}]\boldsymbol{\gamma}_{i}\right)\\
&-\mathrm{vec}\left(\boldsymbol{\gamma}_{i}^{\top}\mathrm{E}[\mathbf{Y}\mathbf{Y}^{\top}]\boldsymbol{\gamma}_{i}\right)\left(\mathrm{E}[\mathbf{Y}]^{\top}\boldsymbol{\gamma}_{i}\right)
+2\left(\boldsymbol{\gamma}_{i}^{\top}\mathrm{E}[\mathbf{Y}]\mathrm{E}[\mathbf{Y}]^{\top}\boldsymbol{\gamma}_{i}\right)\otimes\left(\boldsymbol{\gamma}_{i}^{\top}\mathrm{E}[\mathbf{Y}]\right),
\end{align*}
where $\mathrm{E}[\mathbf{Y}]=\boldsymbol{\mu}+a\boldsymbol{\delta},$  $\mathrm{E}[\mathbf{Y}\mathbf{Y}^{\top}]=\boldsymbol{\mu\mu}^{\top}+a(\boldsymbol{\mu\delta}^{\top}+\boldsymbol{\delta\mu}^{\top})+b\mathbf{\Omega},$ and
\begin{align*}
M_{3}(\mathbf{Y})=&\boldsymbol{\mu}\otimes\boldsymbol{\mu}^{\top}\otimes\boldsymbol{\mu}
 +a[\boldsymbol{\delta}\otimes\boldsymbol{\mu}^{\top}\otimes\boldsymbol{\mu}+\boldsymbol{\mu}\otimes\boldsymbol{\delta}^{\top}\otimes\boldsymbol{\mu}+\boldsymbol{\mu}\otimes\boldsymbol{\mu}^{\top}\otimes\boldsymbol{\delta}]\\
 &+b[\mathbf{\Omega}\otimes\boldsymbol{\mu}+\boldsymbol{\mu}\otimes\mathbf{\Omega}+\mathrm{vec}(\mathbf{\Omega})\otimes\boldsymbol{\mu}^{\top}]\\
&+c[\boldsymbol{\delta}\otimes\mathbf{\Omega}+\mathrm{vec}(\mathbf{\Omega})\boldsymbol{\delta}^{\top}+(\mathbf{I}_{k}\otimes\boldsymbol{\delta})\mathbf{\Omega}-(\mathbf{I}_{k}\otimes\boldsymbol{\delta})(\boldsymbol{\delta}\otimes\boldsymbol{\delta}^{\top})].
\end{align*}

\section{Special cases}
In this section, we consider several special cases within the skew elliptical family, including skew normal, skew $t$, skew logistic, skew Laplace, skew Pearson type II and skew Pearson type VII distributions. Because the measures of skewness and kurtosis for skew elliptical distributions have the same form, we only give $a,~b,~c,~d$, $m_{\ast0}$ and $h_{\ast}$ for different distributions.
\begin{example}\label{ex.1} (Skew normal distribution) Let $\mathbf{Y}\sim SN_{k} (\boldsymbol{\mu},\mathbf{\Omega},\boldsymbol{\delta})$. In this case,
$$g^{(k+1)}(u)=(2\pi)^{-\frac{k+1}{2}}\exp\left(-\frac{u}{2}\right);$$
then
$$\mathrm{E}[R]=\frac{\sqrt{2}\Gamma\left(\frac{k}{2}+1\right)}{\Gamma\left(\frac{k+1}{2}\right)},~\mathrm{E}[R^{2}]=k+1,~\mathrm{E}[R^{3}]=\frac{2\sqrt{2}\Gamma\left(\frac{k}{2}+2\right)}{\Gamma\left(\frac{k+1}{2}\right)},~\mathrm{E}[R^{4}]=(k+1)(k+3).$$
Further,
$$a=c=\sqrt{\frac{2}{\pi}}~ \mathrm{and} ~b=d=1.$$
In addition, $m_{\ast0}$ is the unique solution of the equation
\begin{align*}
y_{1}\Phi(\eta(y_{1}))-\phi(\eta(y_{1}))=0
\end{align*}
and
\begin{align*}
h_{\ast}=-a\delta_{\ast}+\frac{\eta(1)\phi(\eta(a\delta_{\ast}))}{\Phi(\eta(a\delta_{\ast}))},
\end{align*}
where $\Phi(\cdot)$ and $\phi(\cdot)$ are the cdf and pdf of standard normal distribution, respectively, and
$\eta(y_{1})=\frac{\delta_{\ast}y_{1}}{\sqrt{1-\delta_{\ast}^{2}}}$.

\end{example}

Specifically, we can further simplify the measures as
\begin{align*}
\beta_{1,k}
 &=\frac{2(4-\pi)^{2}\delta_{\ast}^{6}}{(\pi-2\delta_{\ast}^{2})^{3}},
\\
 \beta_{2,k}
 &=\frac{\pi(3\pi-12\delta_{\ast}^{2}-4\delta_{\ast}^{4})}{(\pi-2\delta_{\ast}^{2})^{2}}
 +k^{2}+2k-3.
 \end{align*}
 \begin{remark}
 Let $\phi_{k}$ and $\Phi$ denote pdf of $k$-dimensional standard normal distribution and the cumulative distribution function of the univariate standard normal distribution, respectively. Then, Zeller et al. (2011) and Louredo et al. (2022) introduced the skew normal random vector $\mathbf{Y}\sim SN_{k} (\boldsymbol{\mu},\mathbf{\Omega},\boldsymbol{\lambda})$ with its pdf as
 $$f_{\mathbf{Y}}(\boldsymbol{y})=2\phi_{k}(\boldsymbol{y}|\boldsymbol{\mu},\mathbf{\Omega})\Phi(\boldsymbol{\lambda}^{\top}\mathbf{\Omega}^{-1/2}(\boldsymbol{y-\mu})).$$
 We find that there is a relationship between the parameter $\boldsymbol{\lambda}$ above
  and the parameter $\boldsymbol{\delta}$ as follows:
  \begin{align}\label{(a20)}
  \boldsymbol{\delta}^{\top}=\frac{\boldsymbol{\lambda}^{\top}\Omega^{1/2}}{\sqrt{1+\boldsymbol{\lambda}^{\top}\boldsymbol{\lambda}}}.
  \end{align}
\end{remark}
\begin{example}\label{ex.2} (Skew $t$ distribution) Let $\mathbf{Y}\sim St_{k} (\boldsymbol{\mu},\mathbf{\Omega},\boldsymbol{\delta},m)$. In this case,
$$g^{(k+1)}(u)=\frac{\Gamma\left(\frac{m+k+1}{2}\right)}{\Gamma\left(\frac{m}{2}\right)(m\pi)^{(k+1)/2}}\left(1+\frac{u}{m}\right)^{-\frac{m+k+1}{2}};$$
then
\begin{align*}
&\mathrm{E}[R]=\frac{\sqrt{m}\Gamma\left(\frac{k}{2}+1\right)\Gamma\left(\frac{m-1}{2}\right)}{\Gamma\left(\frac{k+1}{2}\right)\Gamma\left(\frac{m}{2}\right)}~(~m>1),
~\mathrm{E}[R^{2}]=\frac{m(k+1)}{m-2}~(m>2),\\
&\mathrm{E}[R^{3}]=\frac{m^{3/2}\Gamma\left(\frac{k}{2}+2\right)\Gamma\left(\frac{m-3}{2}\right)}{\Gamma\left(\frac{k+1}{2}\right)\Gamma\left(\frac{m}{2}\right)}~(m>3),
~\mathrm{E}[R^{4}]=\frac{m^{2}(k+1)(k+3)}{(m-2)(m-4)}~(m>4),
\end{align*}
\begin{align*}
&a=\frac{\sqrt{m}\Gamma\left(\frac{m-1}{2}\right)}{\sqrt{\pi}\Gamma\left(\frac{m}{2}\right)}~(m>1),
~b=\frac{m}{m-2}~(m>2),
\\
&c=\frac{m^{3/2}\Gamma\left(\frac{m-3}{2}\right)}{2\sqrt{\pi}\Gamma\left(\frac{m}{2}\right)}~(m>3),
~d=\frac{m^{2}}{(m-2)(m-4)}~(m>4).
\end{align*}
Furthermore,
$m_{\ast0}$ is the unique solution of the equation
\begin{align*}
-\frac{(m+k+1)y_{1}}{m}\int_{-\infty}^{\eta(y_{1})}
\left(1+\frac{r^{2}+y_{1}^{2}}{m}\right)^{-(m+k+3)/2}\mathrm{d}r+\eta(1)\left(1+\frac{\eta^{2}(y_{1})+y_{1}^{2}}{m}\right)^{-(m+k+1)/2}=0
\end{align*}
and
\begin{align*}
h_{\ast}=\frac{-\frac{a\delta_{\ast}(m+k+1)}{m}\int_{-\infty}^{\eta(a\delta_{\ast})}
\left(1+\frac{r^{2}+a^{2}\delta_{\ast}^{2}}{m}\right)^{-\frac{m+k+3}{2}}\mathrm{d}r
+\eta(1)\left(1+\frac{(\eta(a\delta_{\ast}))^{2}+a^{2}\delta_{\ast}^{2}}{m}\right)^{-\frac{m+k+1}{2}}}{\int_{-\infty}^{\eta(a\delta_{\ast})}
\left(1+\frac{r^{2}+a^{2}\delta_{\ast}^{2}}{m}\right)^{-\frac{m+k+1}{2}}\mathrm{d}r},
\end{align*}
with $\eta(y_{1})=\frac{\delta_{\ast}y_{1}}{\sqrt{1-\delta_{\ast}^{2}}}$.
\end{example}

\begin{example}\label{ex.3} (Skew logistic distribution) Let $\mathbf{Y}\sim SLo_{k} (\boldsymbol{\mu},\mathbf{\Omega},\boldsymbol{\delta})$. In this case,
$$g^{(k+1)}(u)=\frac{\exp(-u)}{1+\exp(-u)};$$
then \begin{align*}
&\mathrm{E}[R]=\frac{\pi^{(k+1)/2}\Gamma\left(\frac{k}{2}+1\right)}{\Gamma\left(\frac{k+1}{2}\right)}\Psi_{1}^{\ast}\left(-1,\frac{k}{2}+1,1\right),~
\mathrm{E}[R^{2}]=\frac{\pi^{(k+1)/2}(k+1)}{2}\Psi_{1}^{\ast}\left(-1,\frac{k+3}{2},1\right),\\
&\mathrm{E}[R^{3}]=\frac{\pi^{(k+1)/2}\Gamma\left(\frac{k}{2}+2\right)}{\Gamma\left(\frac{k+1}{2}\right)}\Psi_{1}^{\ast}\left(-1,\frac{k}{2}+2,1\right),~
\mathrm{E}[R^{4}]=\frac{\pi^{(k+1)/2}(k+3)(k+1)}{4}\Psi_{1}^{\ast}\left(-1,\frac{k+5}{2},1\right).
\end{align*}
Further,
\begin{align*}
&a=\pi^{k/2}\Psi_{1}^{\ast}\left(-1,\frac{k}{2}+1,1\right),
~b=\frac{\pi^{(k+1)/2}}{2}\Psi_{1}^{\ast}\left(-1,\frac{k+3}{2},1\right),
\\
&c=\frac{\pi^{k/2}}{2}\Psi_{1}^{\ast}\left(-1,\frac{k}{2}+2,1\right),
~d=\frac{\pi^{(k+1)/2}}{4}\Psi_{1}^{\ast}\left(-1,\frac{k+5}{2},1\right),
\end{align*}
where $\Psi_{\tau}^{\ast}(\varsigma,s,z)$  is the generalized Hurwitz-Lerch zeta function (see Lin et al., 2006).
Moreover, $m_{\ast0}$ is the unique solution of the equation
\begin{align*}
-2y_{1}\int_{-\infty}^{\eta(y_{1})}
\frac{\exp\left\{-\left(r^{2}+y_{1}^{2}\right)\right\}}{\left[1+\exp\left\{-\left(r^{2}+y_{1}^{2}\right)\right\}\right]^{2}}\mathrm{d}r+\eta(1)\frac{\exp\left\{-\left(\eta^{2}(y_{1})+y_{1}^{2}\right)\right\}}{1+\exp\left\{-\left(\eta^{2}(y_{1})+y_{1}^{2}\right)\right\}}=0
\end{align*}
with $\eta(y_{1})=\frac{\delta_{\ast}y_{1}}{\sqrt{1-\delta_{\ast}^{2}}}$, and
\begin{align*}
h_{\ast}=\frac{-2a\delta_{\ast}\int_{-\infty}^{\eta(a\delta_{\ast})}
\frac{\exp\left(-r^{2}-a^{2}\delta_{\ast}^{2}\right)}{\left[1+\exp\left(-r^{2}-a^{2}\delta_{\ast}^{2}\right)\right]^{2}}\mathrm{d}r
+\eta(1)\frac{\exp\left(-(\eta(a\delta_{\ast}))^{2}-a^{2}\delta_{\ast}^{2}\right)}{1+\exp\left(-(\eta(a\delta_{\ast}))^{2}-a^{2}\delta_{\ast}^{2}\right)}}{\int_{-\infty}^{\eta(a\delta_{\ast})}
\frac{\exp\left(-r^{2}-a^{2}\delta_{\ast}^{2}\right)}{1+\exp\left(-r^{2}-a^{2}\delta_{\ast}^{2}\right)}\mathrm{d}r}.
\end{align*}
\end{example}

\begin{example}\label{ex.3} (Skew Laplace distribution) Let $\mathbf{Y}\sim SLa_{k} (\boldsymbol{\mu},\mathbf{\Omega},\boldsymbol{\delta})$. In this case,
$$g^{(k+1)}(u)=\frac{\Gamma\left(\frac{k+1}{2}\right)}{2\pi^{(k+1)/2}\Gamma(k+1)}\exp(-u^{1/2}),~u>0;$$
then
\begin{align*}
&\mathrm{E}[R]=k+1,~
\mathrm{E}[R^{2}]=(k+1)(k+2),\\
&\mathrm{E}[R^{3}]=(k+1)(k+2)(k+3),~
\mathrm{E}[R^{4}]=(k+1)(k+2)(k+3)(k+4).
\end{align*}
Further,
\begin{align*}
&a=\frac{(k+1)\Gamma\left(\frac{k+1}{2}\right)}{\sqrt{\pi}\Gamma\left(\frac{k}{2}+1\right)},
~b=k+2,
\\
&c=\frac{(k+1)(k+2)(k+3)\Gamma\left(\frac{k+1}{2}\right)}{2\sqrt{\pi}\Gamma\left(\frac{k}{2}+2\right)},
~d=(k+2)(k+4).
\end{align*}
Moreover, $m_{\ast0}$ is the unique solution of the equation
\begin{align*}
-y_{1}\int_{-\infty}^{\eta(y_{1})}
\left(r^{2}+y_{1}^{2}\right)^{-1/2}\exp\left\{-\left(r^{2}+y_{1}^{2}\right)^{1/2}\right\}\mathrm{d}r+\eta(1)\exp\left\{-\left(\eta^{2}(y_{1})+y_{1}^{2}\right)^{1/2}\right\}=0
\end{align*}
and
\begin{align*}
h_{\ast}=\frac{-a\delta_{\ast}\int_{-\infty}^{\eta(a\delta_{\ast})}
\left(r^{2}+a^{2}\delta_{\ast}^{2}\right)^{-1/2}\exp\left\{-\left(r^{2}+a^{2}\delta_{\ast}^{2}\right)^{1/2}\right\}\mathrm{d}r
+\eta(1)\exp\left\{-\left((\eta(a\delta_{\ast}))^{2}+a^{2}\delta_{\ast}^{2}\right)^{1/2}\right\}}{\int_{-\infty}^{\eta(a\delta_{\ast})}
\exp\left\{-\left(r^{2}+a^{2}\delta_{\ast}^{2}\right)^{1/2}\right\}\mathrm{d}r},
\end{align*}
with $\eta(y_{1})=\frac{\delta_{\ast}y_{1}}{\sqrt{1-\delta_{\ast}^{2}}}$.
\end{example}

\begin{example}\label{ex.4} (Skew Pearson type II distribution) Let $\mathbf{Y}\sim SPII_{k} (\boldsymbol{\mu},\mathbf{\Omega},\boldsymbol{\delta},t)$. In this case,
$$g^{(k+1)}(u)=\frac{\Gamma\left(t+1+\frac{k+1}{2}\right)}{\Gamma\left(t+1\right)\pi^{(k+1)/2}}\left(1-u\right)^{t},~0<u<1,~t>-1;$$
 then
 \begin{align*}
&\mathrm{E}[R]=\frac{\Gamma\left(\frac{k+1}{2}+t+1\right)\Gamma\left(\frac{k}{2}+1\right)}{\Gamma\left(\frac{k}{2}+t+2\right)\Gamma\left(\frac{k+1}{2}\right)},
~\mathrm{E}[R^{2}]=\frac{k+1}{k+2t+3},\\
&\mathrm{E}[R^{3}]=\frac{\Gamma\left(\frac{k+1}{2}+t+1\right)\Gamma\left(\frac{k}{2}+2\right)}{\Gamma\left(\frac{k}{2}+t+3\right)\Gamma\left(\frac{k+1}{2}\right)},
~\mathrm{E}[R^{4}]=\frac{\left(k+1\right)\left(k+3\right)}{\left(k+2t+3\right)\left(k+2t+5\right)}.
\end{align*}
Further,
\begin{align*}
&a=\frac{\Gamma\left(t+1+\frac{k+1}{2}\right)}{\sqrt{\pi}\Gamma\left(t+2+\frac{k}{2}\right)},
~b=\frac{1}{2t+k+3},
\\
&c=\frac{\Gamma\left(t+1+\frac{k+1}{2}\right)}{2\sqrt{\pi}\Gamma\left(t+3+\frac{k}{2}\right)},
~
d=\frac{1}{(2t+k+3)(2t+k+5)}.
\end{align*}
Moreover, $m_{\ast0}$ is the unique solution of the equation
\begin{align*}
-2ty_{1}\int_{-\infty}^{\eta(y_{1})}
\left(1-r^{2}-y_{1}^{2}\right)^{t-1}\mathrm{d}r+\eta(1)g^{(k+1)}\left(1-\eta^{2}(y_{1})-y_{1}^{2}\right)^{t}=0
\end{align*}
and
\begin{align*}
h_{\ast}=\frac{-2at\delta_{\ast}\int_{-\infty}^{\eta(a\delta_{\ast})}
\left(1-(r^{2}+a^{2}\delta_{\ast}^{2})\right)^{t-1}\mathrm{d}r
+\eta(1)\left(1-((\eta(a\delta_{\ast}))^{2}+a^{2}\delta_{\ast}^{2})\right)^{t}}{\int_{-\infty}^{\eta(a\delta_{\ast})}
\left(1-(r^{2}+a^{2}\delta_{\ast}^{2})\right)^{t}\mathrm{d}r},
\end{align*}
with $\eta(y_{1})=\frac{\delta_{\ast}y_{1}}{\sqrt{1-\delta_{\ast}^{2}}}$.
\end{example}

\begin{example}\label{ex.5} (Skew Pearson type VII distribution) Let $\mathbf{Y}\sim SPVII_{k} (\boldsymbol{\mu},\mathbf{\Omega},\boldsymbol{\delta},t)$. In this case,
$$g^{(k+1)}(u)=\frac{\Gamma(t)}{\Gamma\left(t-\frac{k+1}{2}\right)\pi^{(k+1)/2}}\left(1+u\right)^{-t},~t>\frac{k+1}{2};$$
 then
  \begin{align*}
&\mathrm{E}[R]=\frac{\Gamma\left(\frac{k}{2}+1\right)\Gamma\left(t-\frac{k}{2}-1\right)}{\Gamma\left(\frac{k+1}{2}\right)\Gamma\left(t-\frac{k+1}{2}\right)},
~\mathrm{E}[R^{2}]=\frac{k+1}{2t-k-3},\\
&\mathrm{E}[R^{3}]=\frac{\Gamma\left(\frac{k}{2}+2\right)\Gamma\left(t-\frac{k}{2}-2\right)}{\Gamma\left(\frac{k+1}{2}\right)\Gamma\left(t-\frac{k+1}{2}\right)},
~\mathrm{E}[R^{4}]=\frac{\left(k+1\right)\left(k+3\right)}{\left(2t-k-3\right)\left(2t-k-5\right)};
\end{align*}
Further,
\begin{align*}
&a=\frac{\Gamma\left(t-\frac{k}{2}-1\right)}{\sqrt{\pi}\Gamma\left(t-\frac{k+1}{2}\right)}~\left(t>\frac{k}{2}+1\right),
~
b=\frac{1}{2t-k-3}~\left(t>\frac{k+3}{2}\right),
\\
&c=\frac{\Gamma\left(t-\frac{k}{2}-2\right)}{2\sqrt{\pi}\Gamma\left(t-\frac{k+1}{2}\right)}~\left(t>\frac{k}{2}+2\right),
~
d=\frac{1}{(2t-k-3)(2t-k-5)}~\left(t>\frac{k+5}{2}\right).
\end{align*}
Moreover, $m_{\ast0}$ is the unique solution of the equation
\begin{align*}
-2ty_{1}\int_{-\infty}^{\eta(y_{1})}
\left(1+r^{2}+y_{1}^{2}\right)^{-(t+1)}\mathrm{d}r+\eta(1)\left(1+\eta^{2}(y_{1})+y_{1}^{2}\right)^{-t}=0
\end{align*}
and
\begin{align*}
h_{\ast}=\frac{-2at\delta_{\ast}\int_{-\infty}^{\eta(a\delta_{\ast})}
\left(1+r^{2}+a^{2}\delta_{\ast}^{2}\right)^{-(t+1)}\mathrm{d}r
+\eta(1)\left(1+(\eta(a\delta_{\ast}))^{2}+a^{2}\delta_{\ast}^{2}\right)^{-t}}{\int_{-\infty}^{\eta(a\delta_{\ast})}
\left(1+r^{2}+a^{2}\delta_{\ast}^{2}\right)^{-t}\mathrm{d}r},
\end{align*}
with $\eta(y_{1})=\frac{\delta_{\ast}y_{1}}{\sqrt{1-\delta_{\ast}^{2}}}$.
\end{example}
\section{Application to hypothesis teating}
To construct test statistics based on a given sample, we need to the following lemma.
\begin{lemma}\label{le.2}  Let $\mathbf{Y}\sim SE_{k} (\boldsymbol{\mu},\mathbf{\Omega},\boldsymbol{\delta}, g^{(k+1)})$,  and $\mathbf{C}$ be a $k\times 1$ vector. Then,
\begin{align}\label{(a21)}
\mathbf{C}^{\top}\mathbf{Y}\sim SE_{1} (\mu_{\ast}, \Omega_{\ast},\delta_{\ast}, g^{(2,k+1)}),
\end{align}
where $\mu_{\ast}=\mathbf{C}^{\top}\boldsymbol{\mu}$, $\Omega_{\ast}=\mathbf{C}^{\top}\mathbf{\Omega}\mathbf{C}$, $\delta_{\ast}=\mathbf{C}^{\top}\boldsymbol{\delta}$, and $g^{(2,k+1)}$ is the $2$-dimensional marginal density generator of $g^{(k+1)}$.
\end{lemma}
\noindent $\mathbf{Proof.}$ Let $\lambda=0$ and $m_{i}=1$ in Proposition 4 of Fang (2003). Then we immediately obtain the desired result. $\hfill\square$

As in Malkovich and Afifi (1973), we assume that the test
statistics are computed from a random sample $\mathbf{Y}_{j},~ j = 1, \cdots, n > k$, with the property that every $(k-1)$-dimensional hyperplane has probability zero. We also assume that $\mathbf{C}\neq\boldsymbol{0}$.

 A hypothesis of no skewness will be accepted for the
univariate distribution of $\mathbf{C}^{\top}\mathbf{Y}$ if $b_{1}(\mathbf{C})\leq K _{b_{1}}$,
where $K _{b_{1}}$ is a constant and $$b_{1}(\mathbf{C})=n\frac{\left[\sum_{j=1}^{n}\left(Z_{j}-\overline{Z}\right)^{3}\right]^{2}}{\left[\sum_{j=1}^{n}\left(Z_{j}-\overline{Z}\right)^{2}\right]^{3}},$$
with $Z_{j}=\mathbf{C}^{\top}\mathbf{Y}_{j}$ and $\overline{Z}=\frac{1}{n}\sum_{i=1}^{n}Z_{i}$.

Then, using Roy's principle (see Roy, 1957; SenGupta, 2007), the hypothesis of no multivariate
skewness of the distribution of $\mathbf{Y}$ is accepted if
$$b_{1}^{\ast}=\max_{\mathbf{C}} b_{1}(\mathbf{C})\leq K _{b_{1}}.$$

For kurtosis, we can similarly consider
 a hypothesis of no multivariate kurtosis to be
accepted if
$$[b_{2}^{\ast}]^{2}=\max_{\mathbf{C}} [b_{2}(\mathbf{C})-K]^{2}\leq K _{b_{2}},$$
where
$$b_{2}(\mathbf{C})=n\frac{\sum_{j=1}^{n}\left(Z_{j}-\overline{Z}\right)^{4}}{\left[\sum_{j=1}^{n}\left(Z_{j}-\overline{Z}\right)^{2}\right]^{2}}$$
and $K$ and $K _{b_{ 2}}$ are constants.
We also find that
$$[b_{2}^{\ast}]^{2}=\max\left\{[b_{2}^{\ast}(\max)-K]^2,[b_{2}^{\ast}(\min)-K]^2\right\},$$
where
$$b_{2}^{\ast}(\max)=\max_{\mathbf{C}} b_{2}(\mathbf{C})~\mathrm{and}~b_{2}^{\ast}(\min)=\min_{\mathbf{C}} b_{2}(\mathbf{C}).$$
Note that the statistics $b_{1}^{\ast}$ and $[b_{2}^{\ast}]^{2}$  are invariant with respect to transformations of the form $\mathbf{X}_{j}=\mathbf{A}\mathbf{Y}_{j}+\boldsymbol{b}$, where $\mathbf{A}$ is a nonsingular $k\times k$ matrix and $\boldsymbol{b}$ is any $k\times1$ vector of constants. To simplify computations, we set
$$\mathbf{A}=\sum_{j=1}^{n}(\mathbf{Y}_{j}-\overline{\mathbf{Y}})(\mathbf{Y}_{j}-\overline{\mathbf{Y}})^{\top},$$
where $\overline{\mathbf{Y}}=\frac{1}{n}\sum_{j=1}^{n}\mathbf{Y}_{j}$, and then set $\mathbf{X}_{j}=(\mathbf{A}^{\ast})^{\top}(\mathbf{Y}_{j}-\overline{\mathbf{Y}})$, where $\mathbf{A}^{\ast}$ is such that
$(\mathbf{A}^{\ast})^{\top}\mathbf{A}\mathbf{A}^{\ast}= \mathbf{I}$.

Without loss of
generality, we also restrict the vectors $\mathbf{C}$ to $\mathbf{C}^{\top}\mathbf{C} = 1$. Then, the multivariate sample statistics $b_{1}^{\ast}$ and $[b_{2}^{\ast}]^{2}$ become
$$b_{1}^{\ast}=\max_{\mathbf{C}^{\top}\mathbf{C}=1}n\left[\sum_{j=1}^{n}\left(\mathbf{C}^{\top}\mathbf{X}_{j}\right)^{3}\right]^{2}$$
and
$$[b_{2}^{\ast}]^{2}=\max_{\mathbf{C}^{\top}\mathbf{C}=1}\left[n\sum_{j=1}^{n}\left(\mathbf{C}^{\top}\mathbf{X}_{j}\right)^{4}-K\right]^{2}.$$
The above two maximization problems can be simplified to following two extremal problems:
\begin{align*}
 \begin{cases}
         u(\mathbf{C})=\sum_{j=1}^{n}\left(\mathbf{C}^{\top}\mathbf{X}_{j}\right)^{3},\\
      \mathbf{C}^{\top}\mathbf{C}-1=0,
 \end{cases}
 \end{align*}
 and
 \begin{align*}
 \begin{cases}
         v(\mathbf{C})=\sum_{j=1}^{n}\left(\mathbf{C}^{\top}\mathbf{X}_{j}\right)^{4},\\
      \mathbf{C}^{\top}\mathbf{C}-1=0.
 \end{cases}
 \end{align*}
 Using Lagrange multiplier method, we then have
 \begin{align*}
 L(\mathbf{C},\lambda)=u(\mathbf{C})-\lambda(\mathbf{C}^{\top}\mathbf{C}-1),\\
 \mathcal{L}(\mathbf{C},\gamma)=v(\mathbf{C})-\gamma(\mathbf{C}^{\top}\mathbf{C}-1).
  \end{align*}
  Thus, we get
  \begin{align}\label{(a22)}
  \begin{cases}
 \frac{\partial L(\mathbf{C},\lambda)}{\partial\mathbf{C}}=3\sum_{j=1}^{n}\left(\mathbf{C}^{\top}\mathbf{X}_{j}\right)^{2}\mathbf{X}_{j}-\lambda\mathbf{C}=0,\\
 \frac{\partial L(\mathbf{C},\lambda)}{\partial\lambda}=-(\mathbf{C}^{\top}\mathbf{C}-1)=0,
 \end{cases}
 \Leftrightarrow
 \begin{cases}
 \sum_{j=1}^{n}\left(\mathbf{C}^{\top}\mathbf{X}_{j}\right)^{2}\mathbf{X}_{j}-\lambda^{\ast}\mathbf{C}=0,\\
 \mathbf{C}^{\top}\mathbf{C}-1=0,
 \end{cases}
 \end{align}
 and
 \begin{align}\label{(a23)}
  \begin{cases}
 \frac{\partial \mathcal{L}(\mathbf{C},\gamma)}{\partial\mathbf{C}}=4\sum_{j=1}^{n}\left(\mathbf{C}^{\top}\mathbf{X}_{j}\right)^{3}X_{j}-\gamma\mathbf{C}=0,\\
 \frac{\partial \mathcal{L}(\mathbf{C},\gamma)}{\partial\gamma}=-(\mathbf{C}^{\top}\mathbf{C}-1)=0,
 \end{cases}
 \Leftrightarrow
 \begin{cases}
 \sum_{j=1}^{n}\left(\mathbf{C}^{\top}\mathbf{X}_{j}\right)^{3}X_{j}-\gamma^{\ast}\mathbf{C}=0,\\
 \mathbf{C}^{\top}\mathbf{C}-1=0.
 \end{cases}
 \end{align}
Now, Newton-Raphson iterative method may be used for solving (\ref{(a22)}) and (\ref{(a23)}), and let them be denoted by $\widehat{\lambda}^{\ast}$ and $\widehat{\gamma}^{\ast}$, respectively.
 Then, $b_{1}^{\ast}=n\left(\widehat{\lambda}^{\ast}\right)^{2}$ and $\left\{b_{2}^{\ast}(\max)~\mathrm{or} ~b_{2}^{\ast}(\min)\right\}=n\widehat{\gamma}^{\ast}$.

 For (\ref{(a22)}), the initial vector $\mathbf{C}$ may be chosen
as the vector that maximizes $\lambda^{2}=\left[\sum_{j=1}^{n}\left(\mathbf{C}^{\top}\mathbf{X}_{j}\right)^{3}\right]^2$ after
scanning a lattice of values of $\mathbf{C} = (C_{1}, \cdots, C_{k})^{\top}$ with
$0\leq C_{1} \leq1$.

 For (\ref{(a23)}), the initial scanning of $\mathbf{C}$ determines starting values for finding the maximum and minimum of $\gamma=\sum_{j=1}^{n}\left(\mathbf{C}^{\top}\mathbf{X}_{j}\right)^{4}$.
\section{Empirical study}
Here, we consider the skewness and kurtosis measures for different bivariate skew-elliptical distributions taking parameter $\boldsymbol{\mu}=\boldsymbol{0}$ and  different choices of parameters $\mathbf{\Omega}$ and $\boldsymbol{\delta}$. The obtained results are summarized in Tables from 1-1, 2 to 15-1, 2 (Tables 3-1, 2 to 15-1, 2, are presented in Supplementary materials-table).

\begin{table}[!htbp]\footnotesize%
\centering Table 1-1: Skewness and kurtosis measures for some bivariate skew-elliptical distributions with parameter $\boldsymbol{\mu}=\boldsymbol{0}$ and  some choices of parameters $\mathbf{\Omega}$ and $\boldsymbol{\delta}$.\\
\begin{tabular}{cccccccccccccccccccccccc}
  \hline
  \#& \multicolumn{2}{c}{Parameters} & Distribution &\multicolumn{2}{c}{\underline{ Mardia measures}} &\multicolumn{1}{c}{{ Malkovich-Afifi}}&Isogai\\
              &      &      &              &  $\beta_{1,k}$& $\beta_{2,k}$  &   $\beta_{1}^{\ast} $ & $\boldsymbol{\delta}$\\
  \hline
  1&    $\mathbf{\Omega}$& $\left[\begin{array}{cc}2& 1 \\1& 3\end{array}\right] $  & $SN_{2}(\boldsymbol{\mu},\mathbf{\Omega},\boldsymbol{\delta}) $&9.9624$\times10^{-5}$&7.5698&9.9624$\times10^{-5}$&$\left[\begin{array}{cc}0.2 \\1\end{array}\right] $   \\
   &   $\boldsymbol{\delta}$&$\left[\begin{array}{cc}0.2 \\1\end{array}\right] $ & \\
   \hline
 2&   $\mathbf{\Omega}$& $\left[\begin{array}{cc}2& 1 \\1& 3\end{array}\right] $  & $St_{2}(\boldsymbol{\mu},\mathbf{\Omega},\boldsymbol{\delta},5) $& 0.2153&38.0429
&0.1640& $\left[\begin{array}{cc}0.2 \\1\end{array}\right] $ \\
      &$\boldsymbol{\delta}$&$\left[\begin{array}{cc}0.2 \\1\end{array}\right] $ & \\
  \hline
  3   & $\mathbf{\Omega}$& $\left[\begin{array}{cc}2& 1 \\1& 3\end{array}\right] $  & $SLo_{2}(\boldsymbol{\mu},\mathbf{\Omega},\boldsymbol{\delta}) $&3.9460&11.3448&3.0742& $\left[\begin{array}{cc}0.2 \\1\end{array}\right] $  \\
      &$\boldsymbol{\delta}$&$\left[\begin{array}{cc}0.2 \\1\end{array}\right] $ & \\
   \hline
   4   & $\mathbf{\Omega}$& $\left[\begin{array}{cc}2& 1 \\1& 3\end{array}\right] $  & $SLa_{2}(\boldsymbol{\mu},\mathbf{\Omega},\boldsymbol{\delta}) $&0.0576& 69.9447&0.0442& $\left[\begin{array}{cc}0.2 \\1\end{array}\right] $ \\
     &$\boldsymbol{\delta}$&$\left[\begin{array}{cc}0.2 \\1\end{array}\right] $ & \\
   \hline
   5   & $\mathbf{\Omega}$& $\left[\begin{array}{cc}2& 1 \\1& 3\end{array}\right] $  & $SPII_{2}(\boldsymbol{\mu},\mathbf{\Omega},\boldsymbol{\delta},2) $&0.0088&4.1372&0.0062& $\left[\begin{array}{cc}0.2 \\1\end{array}\right] $  \\
     &$\boldsymbol{\delta}$&$\left[\begin{array}{cc}0.2 \\1\end{array}\right] $ &\\
   \hline
   6  & $\mathbf{\Omega}$& $\left[\begin{array}{cc}2& 1 \\1& 3\end{array}\right] $  & $SPVII_{2}(\boldsymbol{\mu},\mathbf{\Omega},\boldsymbol{\delta},4) $&0.2153&16.0892&0.1640& $\left[\begin{array}{cc}0.2 \\1\end{array}\right] $  \\
      &$\boldsymbol{\delta}$&$\left[\begin{array}{cc}0.2 \\1\end{array}\right] $ & \\
   \hline
\end{tabular}
\end{table}

\begin{table}[!htbp]\tiny
\centering {\footnotesize Table 1-2: Skewness and kurtosis measures for some bivariate skew-elliptical distributions with parameter $\boldsymbol{\mu}=\boldsymbol{0}$ and  some choices of parameters $\mathbf{\Omega}$ and $\boldsymbol{\delta}$.}\\
\begin{tabular}{cccccccccccccccccccccccc}
  \hline
  \#& \multicolumn{2}{c}{Parameters} & Distribution &Song &\multicolumn{2}{c}{\underline{ Balakrishnan-Brito-Quiroz}}&M$\acute{o}$ri-Rohatgi-Sz$\acute{e}$kely& Kollo&Srivastava\\
               &      &      &        &$S(f_{\mathbf{Y}})$       & $T$   &$Q_{\ast}$ & $s(\mathbf{Y})$ &$b(\mathbf{Y})$& $s_{1,k}^{2}$\\
  \hline
  1&    $\mathbf{\Omega}$& $\left[\begin{array}{cc}2& 1 \\1& 3\end{array}\right] $  & $SN_{2}(\boldsymbol{\mu},\mathbf{\Omega},\boldsymbol{\delta}) $&2.5021$\times10^{-5}$& $\left[\begin{array}{cc}-3.0238\times10^{-5} \\ 0.0239\end{array}\right]$& 5.7115$\times10^{-4}$&$\left[\begin{array}{cc}-8.0636\times10^{-5} \\ 0.0637
\end{array}\right]$&$\left[\begin{array}{cc}-8.0432\times10^{-5} \\ 0.0636\end{array}\right]$& 1.7219    \\
   &   $\boldsymbol{\delta}$&$\left[\begin{array}{cc}0.2 \\1\end{array}\right] $ &
\\
   \hline
 2&    $\mathbf{\Omega}$& $\left[\begin{array}{cc}2& 1 \\1& 3\end{array}\right] $  & $St_{2}(\boldsymbol{\mu},\mathbf{\Omega},\boldsymbol{\delta},5) $& 0.0521&$\left[\begin{array}{cc} 3.5695\times10^{-4} \\0.3880 \end{array}\right] $ & 0.1506& $\left[\begin{array}{cc}9.5187 \times10^{-4}\\ 1.0347
\end{array}\right] $&$\left[\begin{array}{cc}  0.4788 \\1.0353\end{array}\right] $&2.8037
\\
   &   $\boldsymbol{\delta}$&$\left[\begin{array}{cc}0.2 \\1\end{array}\right] $ & \\
  \hline
  3&    $\mathbf{\Omega}$& $\left[\begin{array}{cc}2& 1 \\1& 3\end{array}\right] $  & $SLo_{2}(\boldsymbol{\mu},\mathbf{\Omega},\boldsymbol{\delta}) $&1.8931&$\left[\begin{array}{cc} 2.3531 \\-17.9186
\end{array}\right] $&326.6140&  $\left[\begin{array}{cc}6.2748 \\-47.7830\end{array}\right] $ &$\left[\begin{array}{cc}-1.6552\\ -38.0525\end{array}\right] $&24690.5000\\
   &   $\boldsymbol{\delta}$&$\left[\begin{array}{cc}0.2 \\1\end{array}\right] $ & \\
   \hline
   4&    $\mathbf{\Omega}$& $\left[\begin{array}{cc}2& 1 \\1& 3\end{array}\right] $  & $SLa_{2}(\boldsymbol{\mu},\mathbf{\Omega},\boldsymbol{\delta}) $&  0.5609&$\left[\begin{array}{cc}8.8128 \times10^{-5}\\0.2069\end{array}\right] $&0.0428&$\left[\begin{array}{cc}2.3501 \times10^{-4}\\0.5518\end{array}\right] $&$\left[\begin{array}{cc}0.2451 \\0.5520\end{array}\right] $&2.1174
  \\
   &   $\boldsymbol{\delta}$&$\left[\begin{array}{cc}0.2 \\1\end{array}\right] $ & \\
   \hline
   5&    $\mathbf{\Omega}$& $\left[\begin{array}{cc}2& 1 \\1& 3\end{array}\right] $  & $SPII_{2}(\boldsymbol{\mu},\mathbf{\Omega},\boldsymbol{\delta},2) $& 2.0361$\times10^{-67}$&$\left[\begin{array}{cc} 1.2653\times10^{-4} \\-0.0598
 \end{array}\right] $&0.0036&$\left[\begin{array}{cc}3.3741\times10^{-4}\\ -0.1596
\end{array}\right] $&$\left[\begin{array}{cc}  -0.1094 \\-0.1596\end{array}\right] $&1.5687
 \\
   &   $\boldsymbol{\delta}$&$\left[\begin{array}{cc}0.2 \\1\end{array}\right] $ & \\
   \hline
   6&    $\mathbf{\Omega}$& $\left[\begin{array}{cc}2& 1 \\1& 3\end{array}\right] $  & $SPVII_{2}(\boldsymbol{\mu},\mathbf{\Omega},\boldsymbol{\delta},4) $& 0.0393&$\left[\begin{array}{cc}3.5695\times10^{-4}\\ 0.3880\end{array}\right] $&  0.1506
&$\left[\begin{array}{cc} 9.5187\times10^{-4}\\ 1.0347\end{array}\right] $&$\left[\begin{array}{cc}0.4788 \\1.0353
\end{array}\right] $&2.8037
 \\
   &   $\boldsymbol{\delta}$&$\left[\begin{array}{cc}0.2 \\1\end{array}\right] $ & \\
   \hline
\end{tabular}
\end{table}

\begin{table}[!htbp]\footnotesize%
\centering Table 2-1: Skewness and kurtosis measures for some bivariate skew-elliptical distributions with parameter $\boldsymbol{\mu}=\boldsymbol{0}$ and  some choices of parameters $\mathbf{\Omega}$ and $\boldsymbol{\delta}$.\\
\begin{tabular}{cccccccccccccccccccccccc}
  \hline
  \#& \multicolumn{2}{c}{Parameters} & Distribution &\multicolumn{2}{c}{\underline{ Mardia measures}} &\multicolumn{1}{c}{{ Malkovich-Afifi}}&Isogai \\
              &      &      &              &  $\beta_{1,k}$& $\beta_{2,k}$  &   $\beta_{1}^{\ast} $ & $\boldsymbol{\delta}$\\
  \hline
  1&    $\mathbf{\Omega}$& $\left[\begin{array}{cc}2& 1 \\1& 3\end{array}\right] $  & $SN_{2}(\boldsymbol{\mu},\mathbf{\Omega},\boldsymbol{\delta}) $&0&8&0& $\left[\begin{array}{cc}0\\0  \end{array}\right] $\\
   &   $\boldsymbol{\delta}$&$\left[\begin{array}{cc}0 \\0\end{array}\right] $ & \\
   \hline
 2&   $\mathbf{\Omega}$& $\left[\begin{array}{cc}2& 1 \\1& 3\end{array}\right] $  & $St_{2}(\boldsymbol{\mu},\mathbf{\Omega},\boldsymbol{\delta},5) $& 0&40
&0&  $\left[\begin{array}{cc}0\\0  \end{array}\right] $ \\
      &$\boldsymbol{\delta}$&$\left[\begin{array}{cc}0 \\0\end{array}\right] $ & \\
  \hline
  3   & $\mathbf{\Omega}$& $\left[\begin{array}{cc}2& 1 \\1& 3\end{array}\right] $  & $SLo_{2}(\boldsymbol{\mu},\mathbf{\Omega},\boldsymbol{\delta}) $&0&4.9812&0& $\left[\begin{array}{cc}0\\0  \end{array}\right] $  \\
      &$\boldsymbol{\delta}$&$\left[\begin{array}{cc}0 \\0\end{array}\right] $ & \\
   \hline
   4   & $\mathbf{\Omega}$& $\left[\begin{array}{cc}2& 1 \\1& 3\end{array}\right] $  & $SLa_{2}(\boldsymbol{\mu},\mathbf{\Omega},\boldsymbol{\delta}) $&0&
79.5&0&$\left[\begin{array}{cc}0\\0  \end{array}\right] $  \\
     &$\boldsymbol{\delta}$&$\left[\begin{array}{cc}0 \\0\end{array}\right] $ & \\
   \hline
   5   & $\mathbf{\Omega}$& $\left[\begin{array}{cc}2& 1 \\1& 3\end{array}\right] $  & $SPII_{2}(\boldsymbol{\mu},\mathbf{\Omega},\boldsymbol{\delta},2) $&0&4.1212
&0& $\left[\begin{array}{cc}0\\0  \end{array}\right] $  \\
     &$\boldsymbol{\delta}$&$\left[\begin{array}{cc}0\\0\end{array}\right] $ &\\
   \hline
   6  & $\mathbf{\Omega}$& $\left[\begin{array}{cc}2& 1 \\1& 3\end{array}\right] $  & $SPVII_{2}(\boldsymbol{\mu},\mathbf{\Omega},\boldsymbol{\delta},4) $&0&16&0&$\left[\begin{array}{cc}0\\0  \end{array}\right] $   \\
      &$\boldsymbol{\delta}$&$\left[\begin{array}{cc}0 \\0\end{array}\right] $ & \\
   \hline
\end{tabular}
\end{table}

\begin{table}[!htbp]\footnotesize%
\centering Table 2-2: Skewness and kurtosis measures for some bivariate skew-elliptical distributions with parameter $\boldsymbol{\mu}=\boldsymbol{0}$ and  some choices of parameters $\mathbf{\Omega}$ and $\boldsymbol{\delta}$.\\
\begin{tabular}{cccccccccccccccccccccccc}
  \hline
  \#& \multicolumn{2}{c}{Parameters} & Distribution &Song &\multicolumn{2}{c}{\underline{ Balakrishnan-Brito-Quiroz}}&M$\acute{o}$ri-Rohatgi-Sz$\acute{e}$kely& Kollo&Srivastava\\
               &      &      &        &$S(f_{\mathbf{Y}})$       & $T$   &$Q_{\ast}$ & $s(\mathbf{Y})$ &$b(\mathbf{Y})$& $s_{1,k}^{2}$\\
  \hline
  1&    $\mathbf{\Omega}$& $\left[\begin{array}{cc}2& 1 \\1& 3\end{array}\right] $  & $SN_{2}(\boldsymbol{\mu},\mathbf{\Omega},\boldsymbol{\delta}) $&0& $\left[\begin{array}{cc}0 \\0 \end{array}\right]$& 0&$\left[\begin{array}{cc}0 \\0 \end{array}\right]$&$\left[\begin{array}{cc}0\\ 0\end{array}\right]$& 0   \\
   &   $\boldsymbol{\delta}$&$\left[\begin{array}{cc}0 \\0\end{array}\right] $ &
\\
   \hline
 2&    $\mathbf{\Omega}$& $\left[\begin{array}{cc}2& 1 \\1& 3\end{array}\right] $  & $St_{2}(\boldsymbol{\mu},\mathbf{\Omega},\boldsymbol{\delta},5) $& 0&$\left[\begin{array}{cc} 0 \\0 \end{array}\right] $ & 0& $\left[\begin{array}{cc}0\\ 0\end{array}\right] $&$\left[\begin{array}{cc} 0 \\0\end{array}\right] $&0
\\
   &   $\boldsymbol{\delta}$&$\left[\begin{array}{cc}0 \\0\end{array}\right] $ & \\
  \hline
  3&    $\mathbf{\Omega}$& $\left[\begin{array}{cc}2& 1 \\1& 3\end{array}\right] $  & $SLo_{2}(\boldsymbol{\mu},\mathbf{\Omega},\boldsymbol{\delta}) $&0&$\left[\begin{array}{cc} 0\\ 0\end{array}\right] $&0&  $\left[\begin{array}{cc}0 \\0\end{array}\right] $ &$\left[\begin{array}{cc}0\\ 0\end{array}\right] $&0\\
   &   $\boldsymbol{\delta}$&$\left[\begin{array}{cc}0\\0\end{array}\right] $ & \\
   \hline
   4&    $\mathbf{\Omega}$& $\left[\begin{array}{cc}2& 1 \\1& 3\end{array}\right] $  & $SLa_{2}(\boldsymbol{\mu},\mathbf{\Omega},\boldsymbol{\delta}) $&  0&$\left[\begin{array}{cc}0 \\0\end{array}\right] $&0&$\left[\begin{array}{cc}0\\ 0\end{array}\right] $&$\left[\begin{array}{cc}0\\ 0\end{array}\right] $&0
&
  \\
   &   $\boldsymbol{\delta}$&$\left[\begin{array}{cc}0 \\0\end{array}\right] $ & \\
   \hline
   5&    $\mathbf{\Omega}$& $\left[\begin{array}{cc}2& 1 \\1& 3\end{array}\right] $  & $SPII_{2}(\boldsymbol{\mu},\mathbf{\Omega},\boldsymbol{\delta},2) $& 0&$\left[\begin{array}{cc} 0 \\0 \end{array}\right] $&0&$\left[\begin{array}{cc}0\\ 0\end{array}\right] $&$\left[\begin{array}{cc} 0\\0\end{array}\right] $&0
 \\
   &   $\boldsymbol{\delta}$&$\left[\begin{array}{cc}0\\0\end{array}\right] $ & \\
   \hline
   6&    $\mathbf{\Omega}$& $\left[\begin{array}{cc}2& 1 \\1& 3\end{array}\right] $  & $SPVII_{2}(\boldsymbol{\mu},\mathbf{\Omega},\boldsymbol{\delta},4) $& 0&$\left[\begin{array}{cc}0\\ 0\end{array}\right] $& 0
&$\left[\begin{array}{cc} 0\\ 0\end{array}\right] $&$\left[\begin{array}{cc}0 \\0\end{array}\right] $&0
 \\
   &   $\boldsymbol{\delta}$&$\left[\begin{array}{cc}0 \\0\end{array}\right] $ & \\
   \hline
\end{tabular}
\end{table}
As can be seen from Tables 1-1, 2 to 15-1, 2, $\beta_{1,k}$ and $\beta_{1}^{\ast} $ for skew normal distribution  have almost the same value, while for other distributions that appear to be quite different.
The $\beta_{2,k}$ value for Laplace distribution appear to be always greater than the corresponding measures of other distributions. Moreover, that values of those measures for skew Pearson type VII and skew-$t$ distrinutions are the same except for $\beta_{2,k}$ and $S(f_{\mathbf{Y}})$.

The vectorial measures $T$ and $s(\mathbf{Y})$ have very similar results for all distributions in terms of skewness directions. In addition, some differences are remarkable in cases for distributions. For examples, in case (Table 12-1, 2), $b(\mathbf{Y})$ (Kollo measure) and $s(\mathbf{Y})$ have similar results for skew normal distribution, while $b(\mathbf{Y})$ (Kollo measure) and $s(\mathbf{Y})$ have very different results for other distributions.

 It can be seen from Tables 13-1, 2 that  the degrees of freedom, $m$, has a great impact on the values
of the measures. The larger the $m$ is,  the value of those measures is smaller.  From Tables 14-1, 2 to 15-1, 2, we find that
when $t$ is increasing, almost all the measures (or elements of vector) for skew Pearson type II distribution are decreasing. In addition, the values of the measures of skew Pearson type VII distribution  are similar to those of skew Pearson type II distribution in $t$.
\section{Real-life data analysis}

We consider the daily log-returns of the closing price of companies in the Consumer Staples and Energy sectors. Consumer Staples sector include Church \& Dwight (CHD), Clorox Company (CLX) and Coca-Cola Company (KO) companies, and Energy sector consists of Conoco Phillips (COP), Devon Energy (DVN) and  Chevron Corp. (CVX) companies (for details, see Amiri and Balakrishnan, 2022, Table 1). The data are from the workdays during the period 2022-04-01 to 2023-03-31 (For the data, see http://www.nasdaq.com/).

Because SE distributions is a wide family, it is difficult to carry out its parameter evaluation.
We use R Package $skewMLRM$ to fit data sets for different subfamilies, and then
we choose the best fit among them. This package enables us to select the model inside the subfamilies, including the multivariate scale mixtures of normal (SMN), the multivariate scale mixtures of
skew-normal (SMSN) and the multivariate skew scale mixtures of normal (SSMN) classes.

Let $\mathbf{U}=(U_{1},U_{2},U_{3})^{\top}$ be the log-returns for the companies CHD, CLX and KO, and $\mathbf{V}=(V_{1},V_{2},V_{3})^{\top}$ be the log-returns for the companies COP, DVN and  CVX, respectively. Then, we use the $skewMLRM$ package to choose the model in each subclass
and their Akaike information criteria (AIC) and Bayesian information criteria
(BIC) values are presented in Table 16.
\begin{table}[!htbp]
\centering Table 16: Fitted distributions and their information criterion values to the logarithm of closing prices of different sectors.\\
\setlength{\tabcolsep}{0.7mm}{
\begin{tabular}{cccccccccccccccc}
  \hline
  \hline
  {Subclass}&\multicolumn{2}{c}{SMN}$\big|$&\multicolumn{2}{c}{SMSN~and~SSMN }\\
  \hline
  Sample &$\mathbf{U}$&$\mathbf{V}$&$\mathbf{U}$&$\mathbf{V}$\\
  Model& N&N&SN&SN\\
  AIC& -2363.844&-1798.6265&-2417.458&-1851.4617\\
  BIC &-2332.115&-1766.8974&-2375.152&-1809.1562\\
  \hline
  \hline
\end{tabular}}
\end{table}

According to the minimum AIC and BIC, the skew-normal distributions provides the best fit to the data. We thus find
$$\mathbf{U}\sim SN_{3} (\boldsymbol{\mu},\mathbf{\Omega},\boldsymbol{\delta}),~\mathbf{V}\sim SN_{3} (\boldsymbol{\mu}',\mathbf{\Omega}',\boldsymbol{\delta}'), $$
and the maximum likelihood estimates of the corresponding model parameters are
\begin{align*}
\boldsymbol{\mu}=\left(\begin{array}{ccccccccccccccccc}
4.48872\\
5.02663\\
4.15431
\end{array}
\right)
,
\mathbf{\Omega}=\left(\begin{array}{ccccccccccccccccccc}
0.10117&   0.02285&    0.02948\\
 0.02285&0.06727&   0.01941\\
0.02948&0.01941 &    0.04078
\end{array}
\right),~
\boldsymbol{\delta}=\left(\begin{array}{ccccccccccccccc}
-0.1494547\\ -0.2270945 \\-0.1434309
\end{array}
\right),
\end{align*}
\begin{align*}
\boldsymbol{\mu}'=\left(\begin{array}{ccccccccccccccccc}
4.71430\\
4.11939\\
5.16496
\end{array}
\right)
,
\mathbf{\Omega}'=\left(\begin{array}{ccccccccccccccccccc}
0.10042&   0.03895&    0.04782\\
 0.03895&0.11667&   0.01062\\
0.04782&0.01062 &    0.08271
\end{array}
\right),~
\boldsymbol{\delta}'=\left(\begin{array}{ccccccccccccccc}
-0.09351741 \\0.09451593 \\-0.26141400
\end{array}
\right).
\end{align*}
PP-plots and bivariate scatter plots, with contour lines, for the fitted (SN) distributions  are presented in Figures 1 and 2.

 Next, we compute different measures of
skewness and  kurtosis for the two sectors (or SN distributions), and the obtained results are displayed in Tables 17-1, 2.
\begin{table}[!htbp]\footnotesize%
\centering Table 17-1: Skewness and kurtosis measures for the two sectors.\\
\begin{tabular}{cccccccccccccccccccccccc}
  \hline
  \#& Distribution &\multicolumn{2}{c}{\underline{ Mardia measures}} &\multicolumn{1}{c}{{ Malkovich-Afifi}}&Isogai\\
                       &              &  $\beta_{1,k}$& $\beta_{2,k}$  &   $\beta_{1}^{\ast} $ & $\boldsymbol{\delta}$\\
  \hline
  1&      $\mathbf{U}\sim SN_{3}(\boldsymbol{\mu},\mathbf{\Omega},\boldsymbol{\delta}) $&0.42214
&17.29465&0.42214&$\left[\begin{array}{c} -0.14945\\ -0.22709 \\-0.14343\end{array}\right] $   \\
   \hline
 2& $\mathbf{V}\sim SN_{3}(\boldsymbol{\mu}',\mathbf{\Omega}',\boldsymbol{\delta}') $&0.66939 &19.12467
&0.66939& $\left[\begin{array}{c}-0.09352 \\0.09452 \\-0.26141\end{array}\right] $ \\
  \hline
\end{tabular}
\end{table}

\begin{table}[!htbp]\footnotesize
\centering Table 17-2: Skewness and kurtosis measures for the two sectors.\\
\begin{tabular}{cccccccccccccccccccccccc}
  \hline
  \# & Distribution &Song &\multicolumn{2}{c}{\underline{ Balakrishnan-Brito-Quiroz}}&M$\acute{o}$ri-Rohatgi-Sz$\acute{e}$kely& Kollo&Srivastava\\
                          &        &$S(f_{\mathbf{Y}})$       & $T$   &$Q_{\ast}$ & $s(\mathbf{Y})$ &$b(\mathbf{Y})$& $s_{1,k}^{2}$\\
  \hline
  1&     $\mathbf{U}\sim SN_{3}(\boldsymbol{\mu},\mathbf{\Omega},\boldsymbol{\delta}) $& 0.17980
& $\left[\begin{array}{c} -0.04137 \\-0.12790 \\-0.08561 \end{array}\right]$& 0.02540 &$\left[\begin{array}{c}-0.20686 \\-0.63950 \\-0.42804
\end{array}\right]$&$\left[\begin{array}{c} -0.52911\\ -1.63568 \\-1.09482
 \end{array}\right]$& 9987.47500   \\
   \hline
 2&     $\mathbf{V}\sim SN_{3}(\boldsymbol{\mu}',\mathbf{\Omega}',\boldsymbol{\delta}') $&  0.23603&$\left[\begin{array}{c} -0.01384 \\ 0.05874 \\-0.16971 \end{array}\right] $ &0.03244
 & $\left[\begin{array}{c}  -0.06919 \\ 0.29371\\ -0.84854
\end{array}\right] $&$\left[\begin{array}{c} -0.03322\\ 0.14101 \\-0.40738\end{array}\right] $&5248.02100
\\
  \hline
\end{tabular}
\end{table}
From Tables 17-1, 2, we observe that $\beta_{1,k}$, $\beta_{2,k}$, $\beta_{1}^{\ast}$, $S(f_{\mathbf{Y}})$, $Q_{\ast}$ and all elements of $b(\mathbf{Y})$,  and the first two elements of $\boldsymbol{\delta}$, $T$ and $s(\mathbf{Y})$  for $\mathbf{V}$ are greater than the corresponding measures for $\mathbf{U}$, while $s_{1,k}^{2}$ and the last elements of $\boldsymbol{\delta}$, $T$ and $s(\mathbf{Y})$  for $\mathbf{V}$ are smaller than the corresponding measures for $\mathbf{U}$.

 \section{Concluding remarks}

 In this paper, we have derived and compared eight measures of
skewness and Mardia measure of kurtosis for skew-elliptical distributions, which is an extension on the work on skew normal distribution by Balakrishnan and Scarpa (2012). We have also derived exact expressions of these measures for skew-elliptical family, and have presented the results for some special cases, such as skew normal, skew $t$, skew logistic, skew Laplace, skew Pearson type II and skew Pearson type VII distributions.
 Recently, Loperfido (2020) investigated some properties of Koziol's measure of multivariate kurtosis, with a motivation for its use in statistical practice.  Chowdhury et al. (2022) proposed sub-dimensional Mardia measures of multivariate skewness and kurtosis, which are the Mardia measures of multivariate skewness and kurtosis for all marginal distributions. We may develop similar results for skew-elliptical distributions. It will also be of interest to develop other kurtosis measures for the family of skew-elliptical distributions. Work on these problems is currently under progress and we hope to report the corresponding findings in a future paper.
\section*{Acknowledgments}

\noindent  The research was supported by the National Natural Science Foundation of China (No. 12071251)
\section*{Conflicts of Interest}
\noindent The authors declare that they have no conflicts of interest.

\medskip

\end{document}